\documentclass[11pt]{article}

\usepackage[tbtags]{amsmath}

\newcommand{\interligne}[1]{\renewcommand{\baselinestretch}{#1} \small\normalsize}

\usepackage{amsthm,amsopn,amsfonts,amssymb,epsfig,pictex}

\edef\savecatcodeat{\the\catcode`@}
\catcode`\@=11

\def\tb@ifSpecChars#1#2{#1}
\def\tb@ifNoSpecChars#1#2{#2}

\def\tableau{%
  \bgroup
  \@ifstar{\let\Tif\tb@ifNoSpecChars\tb@tableauB}
          {\let\Tif\tb@ifSpecChars\tb@tableauB}}

\def\tb@tableauB{
  \@ifnextchar[{\tb@tableauC}{\tb@tableauC[]}}

\def\tb@tableauC[#1]{\hbox\bgroup%
    \let\\=\cr
    \def\bl{\global\let\tbcellF\tb@cellNF}%
    \def\tf{\global\let\tbcellF\tb@cellH}
    \dimen2=\ht\strutbox \advance\dimen2 by\dp\strutbox%
    \ifx\baselinestretch\undefined\relax%
    \else%
       \dimen0=100sp \dimen0=\baselinestretch\dimen0%
       \dimen2=100\dimen2 \divide\dimen2 by\dimen0%
    \fi%
    \let\tpos\tb@vcenter
    \tb@initYoung
    \tb@options#1\eoo
    \let\arrow\tb@arrow%
    \dimen0=\Tscale\dimen2%
    \dimen1=\dimen0 \advance\dimen1 by \tb@fframe%
    \lineskip=0pt\baselineskip=0pt
    \def\tb@nothing{}%
    \def\endcellno{$\rss\egroup\bss\egroup}
    \def\endcell{\endcellno\kern-\dimen0}
    \def\begincell{\vbox to\dimen0\bgroup\vss\hbox to\dimen0\bgroup\hss$}%
    \let\overlay\tb@overlay%
    \let\fl\tb@fl%
    \let\lss\hss\let\rss\hss\let\tss\vss\let\bss\vss
    \def\mkcell##1{
        \let\tbcellF\tb@cellD
        \def\tb@cellarg{##1}
        \ifx\tb@cellarg\tb@nothing\let\tb@cellarg\tb@cellE\fi%
        \begincell\tb@cellarg\endcellno
        \tbcellF}
    \let\savecellF\tbcellF
     \Tif{\catcode`,=4\catcode`|=\active}{}\tb@tableauD}%

\let\tb@savetableauD\tableauD
{
    \catcode`|=\active \catcode`*=\active \catcode`~=\active%
    \catcode`@=\active
\gdef\tableauD#1{%
  \Tif{
    \mathcode`|="8000 \mathcode`*="8000%
    \mathcode`~="8000 \mathcode`@="8000%
    \def@{\bullet}%
    \let|\cr
    \let*\tf
    \let~\sk
  }{}%
  \tpos{\tabskip=0pt\halign{&\mkcell{##}\cr#1\crcr}}%
  \global\let\tbcellF\savecellF
  \egroup
  \egroup}
}
\let\tb@tableauD\tableauD
\let\tableauD\tb@savetableauD
\let\tb@savetableauD\undefined

\def\tb@options#1{\ifx#1\eoo\relax\else\tb@option#1\expandafter\tb@options\fi}

\def\tb@option#1{%
  \if#1t\let\tpos\tb@vtop\fi
  \if#1c\let\tpos\tb@vcenter\fi
  \if#1b\let\tpos\vbox\fi
  \if#1F\tb@initFerrers\fi
  \if#1Y\tb@initYoung\fi
  \if#1s\tb@initSmall\fi
  \if#1m\tb@initMedium\fi
  \if#1l\tb@initLarge\fi
  \if#1p\tb@initPartition\fi
  \if#1a\tb@initArrow\fi
}

\def\tb@vcenter#1{\ifmmode\vcenter{#1}\else$\vcenter{#1}$\fi}

\def\tb@vtop#1{\hbox{\raise\ht\strutbox\hbox{\lower\dimen0\vtop{#1}}}}

\def\tb@initPartition{\def\Tscale{.3}}
\def\tb@initSmall{\def\Tscale{1}}
\def\tb@initMedium{\def\Tscale{2}}
\def\tb@initLarge{\def\Tscale{3}}

\def\tb@initArrow{\dimen2=1.25em}

\def\tb@initYoung{%
  \def\tb@cellE{}
  \let\tb@cellD\tb@cellN
  \def\sk{\global\let\tbcellF\tb@cellNF}}
\def\tb@initFerrers{%
  \def\tb@cellE{\bullet}
  \let\tb@cellD\tb@cellNF
  \def\sk{\bullet}}

\tb@initMedium

\def\tb@sframe#1{%
  \vbox to0pt{
    \vss
    \hbox to0pt{%
      \hss
      \vbox to\dimen1{
        \hrule depth #1 height0pt
        \vss
        \hbox to\dimen1{
          \vrule width #1 height\dimen1
          \hss
          \vrule width #1
          }%
        \vss
        \hrule height #1 depth 0in
        }%
      \kern-\tb@hframe
      }%
    \kern-\tb@hframe}}

\def\tb@hframe{.2pt}\def\tb@fframe{.4pt}\def\tb@bframe{2pt}
\def\tb@cellH{\tb@sframe{\tb@bframe}}       
\def\tb@cellNF{}                            
\def\tb@cellN{\tb@sframe{\tb@fframe}}       
\let\tbcellF\tb@cellN                       

\def\tb@rpad{1pt}
\def\tb@lpad{1pt}
\def\tb@tpad{1.8pt}
\def\tb@bpad{1.8pt}

\def\tb@overlay{\endcell\@ifnextchar[{\tb@overlaya}{\begincell}}
\def\tb@overlaya[#1]{\vbox to\dimen0\bgroup%
  \tb@overlayoptions#1\eoo%
  \tss\hbox to\dimen0\bgroup\lss$}
\def\tb@overlayoptions#1{\ifx#1\eoo\relax\else\tb@overlayoption#1\expandafter\tb@overlayoptions\fi}

\def\tb@overlayoption#1{
  \if#1t\def\tss{\vskip\tb@tpad}\let\bss\vss\fi
  \if#1c\let\tss\vss\let\bss\vss\fi
  \if#1b\def\bss{\vskip\tb@bpad}\let\tss\vss\fi
  \if#1l\def\lss{\hskip\tb@lpad}\let\rss\hss\fi
  \if#1m\let\lss\hss\let\rss\hss\fi
  \if#1r\def\rss{\hskip\tb@rpad}\let\lss\hss\fi
}

\def\tb@fl{\endcell\begincell\vrule depth 0pt width \dimen0 height \dimen0 \endcell\begincell}

\@ifundefined{diagram}{}{
\def\tb@arrowpad{.5}

\newoptcommand{\tb@arrow}{\@ne}[2]{%
  \endcell
   \begingroup%
   \let\dg@getnodesize\tb@getnodesize
   \dg@USERSIZE=#1\relax%
   \ifnum\dg@USERSIZE<\@ne \dg@USERSIZE=\@ne \fi%
   \dg@parse{#2}%
   \dg@label{\tb@draw{#1}{#2}}}

\def\tb@getnodesize#1#2#3#4#5{\dimen3=\tb@arrowpad\dimen2 #4=\dimen3 #5=\dimen3\relax}
\def\tb@getnodesize#1#2#3#4#5{\ifnum#2=0\ifnum#3=0\tb@getnodesizetail{#4}{#5}\else\tb@getnodesizehead{#4}{#5}\fi\else\tb@getnodesizehead{#4}{#5}\fi}
\def\tb@getnodesizetail#1#2{\dimen3=.5\dimen2 #1=\dimen3 #2=\dimen3}
\def\tb@getnodesizehead#1#2{\dimen3=.5\dimen2 #1=\dimen3 #2=\dimen3}

\def\tb@draw#1#2#3#4{%
        \dg@X=0\dg@Y=0\dg@XGRID=1\dg@YGRID=1\unitlength=.001\dimen0%
        \dg@LBLOFF=\dgLABELOFFSET \divide\dg@LBLOFF\unitlength%
        \dg@drawcalc
        \begincell
        \let\lams@arrow\tb@lams@arrow
        \begin{picture}(0,0)\begingroup\dg@draw{#1}{#2}{#3}{#4}\end{picture}%
        \endcell
        \endgroup
        \begincell}
}

\def\tb@lams@arrow#1#2{%
 \lams@firstx\z@\lams@firsty\z@
 \lams@lastx#1\relax\lams@lasty#2\relax
 \lams@center\z@
 \N@false\E@false\H@false\V@false
 \ifdim\lams@lastx>\z@\E@true\fi
 \ifdim\lams@lastx=\z@\V@true\fi
 \ifdim\lams@lasty>\z@\N@true\fi
 \ifdim\lams@lasty=\z@\H@true\fi
 \NESW@false
 \ifN@\ifE@\NESW@true\fi\else\ifE@\else\NESW@true\fi\fi
 \ifH@\else\ifV@\else
  \lams@slope
  \ifnum\lams@tani>\lams@tanii
   \lams@ht\ten@\p@\lams@wd\ten@\p@
   \multiply\lams@wd\lams@tanii\divide\lams@wd\lams@tani
  \else
   \lams@wd\ten@\p@\lams@ht\ten@\p@
   \divide\lams@ht\lams@tanii\multiply\lams@ht\lams@tani
  \fi
 \fi\fi
 \ifH@  \lams@harrow
 \else\ifV@ \lams@varrow
 \else \lams@darrow
 \fi\fi
}

\catcode`\@=\savecatcodeat
\let\savecatcodeat\undefined

\numberwithin{equation}{section}

\newtheorem{theorem}{Theorem}

\newtheorem{proposition}[theorem]{Proposition}
\newtheorem{conjecture}[theorem]{Conjecture}

\newtheorem{corollary}[theorem]{Corollary} 

\newtheorem{definition}[theorem]{Definition}

\newtheorem{property}[theorem]{Property}

\theoremstyle{remark}

\newtheorem*{acknow}{\bf Acknowledgments}

\def\charge{ {\rm {charge}}}

\def\shape{ {\rm {shape}}}

\def\sstab  {  T }
\def\stab  { \mathcal T }

\def\endprf {\square}

\begin{document}

\begin{center}
{\Large\bf Tableau atoms and a new}

\medskip

{\Large \bf Macdonald positivity conjecture 
}

\vskip 2pt
\bigskip\bigskip\bigskip
 L. Lapointe\\ 
Department of Mathematics and Statistics, McGill University \\  
Montr\'eal, Qu\'ebec, Canada, H3A 2K6\\
{\small\texttt{lapointe@math.mcgill.ca }}\\
\bigskip
A. Lascoux \\
Institut Gaspard Monge,
Universit\'e de Marne-la-Vall\'ee\\
 5 Bd Descartes, Champs sur Marne \\
77454 Marne La Vall\'ee, Cedex, France\\
{\small\texttt{Alain.Lascoux@univ-mlv.fr }} \\
\bigskip
 J. Morse \\
Department of Mathematics,
University of Pennsylvania \\
 209 South 33rd Street\\
 Philadelphia, PA 19104, USA \\
{\small\texttt{morsej@math.upenn.edu}} \\
\bigskip

\bigskip

\bigskip


\bigskip\bigskip

\vskip 2pt

\medskip

\begin{abstract}
Let $\Lambda$ be the space of symmetric functions and
$V_k$ be the subspace spanned by the modified Schur functions
$\{S_\lambda[X/(1-t)]\}_{\lambda_1\leq k}$.
We introduce a new family of 
symmetric polynomials, $\{A_{\lambda}^{(k)}[X;t]\}_{\lambda_1\leq k}$, 
constructed from sums of tableaux using the charge statistic.
We conjecture that the polynomials
$A_{\lambda}^{(k)}[X;t]$ form a basis for $V_k$
and that the Macdonald polynomials indexed by partitions 
whose first part is not larger than $k$ expand positively in terms of 
our polynomials.
A proof of this conjecture would not only imply the Macdonald
positivity conjecture, but would substantially
refine it.   Our construction of the $A_\lambda^{(k)}[X;t]$
relies on the use of tableaux combinatorics
and yields various properties and conjectures
on the nature of these polynomials.
Another important development following from our investigation
is that the $A_{\lambda}^{(k)}[X;t]$ seem to
play the same role for $V_k$ as the Schur functions 
do for $\Lambda$.
In particular, this has led us to the discovery of many
generalizations of properties held by the Schur functions, such as 
Pieri and Littlewood-Richardson type coefficients.
 \end{abstract}
\end{center}
\vskip 2pt
\bigskip\bigskip
\bigskip\bigskip\bigskip\bigskip
\vskip 2pt
\eject

\bigskip
\bigskip
\bigskip

\section{Introduction}
\pagestyle{myheadings} \markboth{\hfill{\sc Tableau atoms and a new Macdonald
positivity conjecture}}
{{\sc L. Lapointe, A. Lascoux and J. Morse}\hfill}

We work with the algebra $\Lambda$ of symmetric functions
in the formal alphabet $x_1,x_2,\ldots$ with coefficients
in $\mathbb Q(q,t)$.  We use $\lambda$-ring notation
in our presentation and refer those unfamiliar with this
device to section 2.
It develops that the filtration of $\Lambda$ given by the spaces 
\begin{equation}
V_k = \left\{S_\lambda[X/(1-t)] \right\}_{\lambda_1\leq k} \, , \quad
\;\text{with}\; k\in \mathbb N \, ,
\end{equation}
provides a natural setting for the study of the 
$q,t$-Kostka coefficients, $K_{\lambda\mu}(q,t)$.  In fact, 
this filtration leads to a family of positivity conjectures 
refining the original Macdonald positivity conjecture,
which now holds 
following the proof \cite{[Ha]} of the $n!$ conjecture \cite{[Ga]}.  
To see how this comes about, we first introduce some notation.

We use a modification of the Macdonald integral forms $J_\mu[X;q,t]$ 
that is obtained by  setting
\begin{equation}
H_\mu[X;q,t] = 
J_\mu[X/(1-t);q,t] =
\sum_\lambda K_{\lambda\mu}(q,t) S_\lambda[X]
\,.
\end{equation}
The integral form $J_\mu[X;q,t]$ at $q=0$
reduces to the Hall-Littlewood polynomial,
\begin{equation*}
J_\mu[X;0,t] = Q_\mu[X;t]
\,.
\end{equation*}
We shall also use a modification of $Q_\mu[X;t]$;
\begin{equation}
H_\mu[X;t] = Q_\mu[X/(1-t);t]
= \sum_{\lambda\geq\mu} K_{\lambda\mu}(t) S_\lambda[X]
\,,
\end{equation}
where $K_{\lambda\mu}(t)$ is the Kostka-Foulkes polynomial.

This given, we should note that bases for $V_k$ also include 
the families \cite{[Ma]}
\begin{equation}
\left\{H_\mu[X;t] \right\}_{\mu_1\leq k}
\quad 
\text{and}
\quad 
\left\{H_\mu[X;q,t] \right\}_{\mu_1\leq k}
\,.
\label{bases}
\end{equation}
Our main contribution is the construction of a new family 
\begin{equation}
\{A_\lambda^{(k)}[X;t]\}_{\lambda_1\leq k} \, ,
\end{equation}
which we conjecture forms a basis for $V_k$
and whose elements, in a sense that can be made precise,
constitute the smallest Schur positive components of $V_k$. 
For this reason, we have chosen to call the
$A_\lambda^{(k)}[X;t]$ the atoms of $V_k$.

We begin by outlining the characterization of our atoms
which may be compared to the combinatorial construction
of the Hall-Littlewood polynomials. 
The formal sum, or the set, of all  
semi-standard tableaux (hereafter called tableaux)
with evaluation $\mu$ will be denoted\footnote{
Double fonts are used to distinguish
sets of tableaux or operators on tableaux from
functions.}
by $\mathbb H_\mu$,
with the convention that $\mathbb H_0$ is the empty tableau.
It was shown in \cite{[LS2]} that
\begin{equation}
H_\mu[X;t] \,=\,\digamma \left(\mathbb H_\mu\right)
\,=\, \sum_{T\in \mathbb H_\mu}
t^{\text{charge}(T)}\,S_{\text{shape}(T)}[X]
\, ,
\label{halldef}
\end{equation}
where $\digamma$ is the functional
\begin{equation}
\digamma( \sstab) = t^{\charge(\sstab)} S_{\shape (\sstab ) }[X] \, .
\label{morphism}
\end{equation}
The formal sum $\mathbb H_\mu$ arises from 
a recursive application of promotion 
operators $\mathbb B_r$ such that 
$\mathbb B_r\mathbb H_\lambda\! = \!\mathbb H_{(r,\lambda)}$:
\begin{equation}
\mathbb H_\mu = 
\mathbb B_{\mu_1} \cdots 
\mathbb B_{\mu_{n-1}} 
\mathbb B_{\mu_n} 
\mathbb H_0
\, .
\label{hallrec}
\end{equation}
The operators $\mathbb B_r$ are tableau analogues of the operators
building recursively the Hall-Littlewood polynomials 
presented in \cite{[Ji],[Ga2]}.

To construct the atoms of $V_k$, we introduce a family of filtering operators, 
$\mathbb P_{\lambda^{\to k}}$, which have the effect of
removing certain elements from the sum of tableaux in \ref{hallrec}
\footnote{The effect of the filtering operator 
is to extract $\lambda$-katabolizable tableaux
(see Section~\ref{subsectab} and \cite{[S1]}).}.
That is, given a $k$-bounded partition $\mu$ (a partition $\mu$ such that $\mu_1 \leq k$), the atom $A^{(k)}_{\mu}[X;t]$ is 
\begin{equation}
A^{(k)}_\mu[X;t] \,=\,\digamma \left(\mathbb A^{(k)}_\mu\right)
\,=\, \sum_{T\in \mathbb A^{(k)}_\mu}
t^{\text{charge}(T)}\,S_{\text{shape}(T)}[X]
\,,
\label{defatoma}
\end{equation}
where 
$\mathbb A^{(k)}_\mu$ is the formal sum of tableaux
obtained from 
\begin{equation}
\mathbb A^{(k)}_\mu = 
\mathbb P_{\mu^{\to k}}
\mathbb B_{\mu_1} 
\cdots
\mathbb P_{(\mu_{n-1},\mu_n)^{\to k}}
\mathbb B_{\mu_{n-1}} 
\mathbb P_{(\mu_n)^{\to k}}
\mathbb B_{\mu_n} \mathbb H_{0}
\, .
\label{defatom1}
\end{equation}
Following from this construction is the
expansion,
\begin{equation}
A^{(k)}_\mu[X;t] = 
S_\mu[X] + \sum_{\lambda>\mu} v^{(k)}_{\lambda\mu}(t) 
S_\lambda[X] \, ,
\quad
\text{with}
\; \;
0\subseteq v_{\lambda\mu}^{(k)}(t)\subseteq K_{\lambda\mu}(t) \, ,
\label{surepos}
\end{equation}
where for two 
polynomials $P,Q\in\mathbb Z[q,t]$, we write
$P\subseteq Q$ to mean  
$Q-P\in \mathbb N[q,t]$.

Originally, the atoms were empirically constructed by the idea that they 
could be characterized by \ref{surepos}
and the two following properties: 

$i)$ for $k$-bounded partitions $\lambda$ and $\mu$,
and any non-zero coefficient $c(t) \in \mathbb N[t]$, 
\begin{equation} 
A_{\lambda}^{(k)}[X;t]-c(t) \, A_{\mu}^{(k)}[X;t] \not = \sum_{\nu} 
v_{\nu}(t) \, S_{\nu}[X] \, , 
\quad\text{where}\quad v_{\nu} \in \mathbb N[t]
\end{equation} 

$ii)$ for any $k$-bounded partition $\mu$,
\begin{equation}
H_\mu[X;t] = 
A_\mu^{(k)}[X;t] +
\sum_{\lambda>\mu} K_{\lambda\mu}^{(k)}(t) \, 
A_\lambda^{(k)}[X;t] \, ,
\quad
\text{with}
\;\;
K_{\lambda\mu}^{(k)}(t)\in \mathbb N[t]\,.
\label{Aprime}
\end{equation}
However, our computer experimentation supported the following 
stronger conjecture, which connects the atoms to Macdonald polynomials indexed by $k$-bounded partitions:
\begin{equation}
H_\mu[X;q,t] = \sum_\lambda K_{\lambda\mu}^{(k)}(q,t)\,  A_\lambda^{(k)}[X;t] \, ,
\label{Bp}
\end{equation}
with
\begin{equation}
0\subseteq K_{\lambda\mu}^{(k)}(q,t)\subseteq K_{\lambda\mu}(q,t)
\,.
\label{Bprime}
\end{equation}
\noindent This has been the primary motivation for the
research that led to this work. 
In particular, given the positive expansion in \ref{surepos}, 
property \ref{Bp} with \ref{Bprime} would not only
prove the Macdonald positivity conjecture, but
would constitute a substantial strengthening of it.

It will transpire that the atoms are a natural 
generalization of the Schur functions.  
In fact, our construction of $A_\lambda^{(k)}[X;t]$ 
yields the property that for large $k$ ($k\geq |\lambda|$), 
\begin{equation}
A_\lambda^{(k)}[X;t]=S_\lambda[X]
\,.
\end{equation}
Thus the atoms of $\Lambda=V_\infty$ are none other than the
Schur functions themselves.  
Moreover, computer exploration has revealed that the 
$A_\lambda^{(k)}[X;t]$ have a variety
of remarkable properties extending and
refining well-known properties of Schur functions.
For example, we have observed generalizations of 
Pieri and Littlewood-Richardson rules, a $k$-analogue of the Young Lattice 
induced by the multiplication action of $e_1$, 
and a $k$-analogue of partition conjugation.
Further, we have noticed that the 
atoms satisfy, on any two alphabets $X$ and $Y$, 
\begin{equation*}
A_\lambda^{(k)}[X+Y;t]=\sum_{|\mu|+|\rho|=|\lambda|}
g_{\mu \rho}^\lambda(t)\,A_{\mu}^{(k)}[X;t]\,A_{\rho}^{(k)}[Y;t] \, ,
\quad\text{where}\;\;
g_{\mu\rho}^\lambda(t)\in \mathbb N[t] \, .
\end{equation*}
The positivity of the coefficients $g_{\mu\rho}^{\lambda}(t)$
appearing here is a natural property of Schur functions not shared by the 
Hall-Littlewood or Macdonald functions.  Finally, the atoms of $V_k$,
when embedded in the atoms of $V_{k'}$ for $k'>k$, seem to decompose positively:
\begin{equation}
A_{\lambda}^{(k)}[X;t]=A_{\lambda}^{(k')}[X;t]+ \sum_{\mu > \lambda} 
v_{\mu \lambda}^{(k \to k')}(t) \, A_{\mu}^{(k')}[X;t] \, ,\quad {\rm{where~}} 
v_{\mu \lambda}^{(k \to k')}(t) \in \mathbb N(t) \, .
\end{equation}

The tableaux combinatorics involved in our construction 
and identity \ref{Aprime}
suggest that the atoms provide a natural structure on
the set of tableaux $\mathbb H_{\mu}$.
For example, we have observed that 
for a $k$-bounded partition $\mu$,
$\mathbb H_\mu$ decomposes 
into disjoint subsets $\mathbb A_T^{(k)}$ indexed by their element of 
minimal charge.   Each of these subsets is characterized 
by the fact that its cyclage-cocyclage poset structure  
is isomorphic to that of
$\mathbb A_{\text{shape}(T)}^{(k)}$, and
since
\begin{equation}
\digamma \left( \mathbb A_T^{(k)} \right)
 = t^{\text{charge}(T)}
A_{\text{shape}(T)}^{(k)}[X;t]
\, ,
\label{digaonat}
\end{equation}
we say that 
$\mathbb A_T^{(k)}$ 
is a copy 
of $\mathbb A_{\text{shape}(T)}^{(k)}$.
Therefore, if ${\cal C}_\mu^{(k)}$ is the
collection of tableaux indexing the 
copies that occur in the decomposition of
$\mathbb H_\mu$, we have
\begin{equation}
\mathbb H_\mu
=\sum_{T\in{\cal C}_\mu^{(k)}}
\mathbb A_T^{(k)}
\,.
\label{HinA}
\end{equation}
Note that the tableaux in 
$\mathbb A^{(k)}_T$ have evaluation $\mu$
while those in $\mathbb A_{\text{shape}(T)}^{(k)}$ 
have, from \ref{defatom1}, evaluation given by the shape of $T$.
Now, \ref{digaonat} and \ref{HinA} imply that
the coefficients $K_{\lambda\mu}^{(k)}(t)$
occurring in \ref{Aprime}
are simply given by the formula
\begin{equation}
K_{\lambda\mu}^{(k)}(t) 
=\sum_{
\begin{subarray}{c}
T \in {\cal C}_\mu^{(k)}\\
\text{shape}(T)=\lambda
\end{subarray}
}
t^{\text{charge}(T)}
\,.
\label{kostkafoulkes1}
\end{equation}

Since the promotion operators $\mathbb B_\ell$ acting on 
$\mathbb A_T^{(k)}$ produce collections of tableaux 
of the same evaluation, we examine their 
decomposition into copies as well.  It appears that
\begin{equation}
\mathbb B_\ell 
\mathbb A_T^{(k)}
= 
\sum_{T' \in {\cal E}_{T,\ell}^{(k)}}
\mathbb A_{T'}^{(k)}\,,
\label{tabpieri}
\end{equation}
where ${\cal E}_{T,\ell}^{(k)}$ is a suitable subcollection
of the tableaux $T'$ of shape $\nu$ such that
$\nu/\text{shape}(T)$ is a horizontal $\ell$-strip.
Therefore, formula \ref{tabpieri} may be considered
a refinement of the classical Pieri rules.
In fact, letting $t=1$ and shape$(T)=\lambda$
in \ref{tabpieri}, we have
\begin{equation}
h_\ell[X]
A_\lambda^{(k)}[X;1]
= 
\sum_{\nu \in  E_{\lambda,\ell}^{(k)}}
A_\nu^{(k)}[X;1]
\,,
\label{pieri}
\end{equation}
where $E_{\lambda,\ell}^{(k)}$
is a subset of the collection of shapes $\nu$ such that
$\nu/\lambda$ is a horizontal $\ell$ strip.  
We shall give a simple combinatorial procedure 
for determining $E_{\lambda,\ell}^{(k)}$.

\begin{figure}[htb]
\begin{center}
\epsfig{file=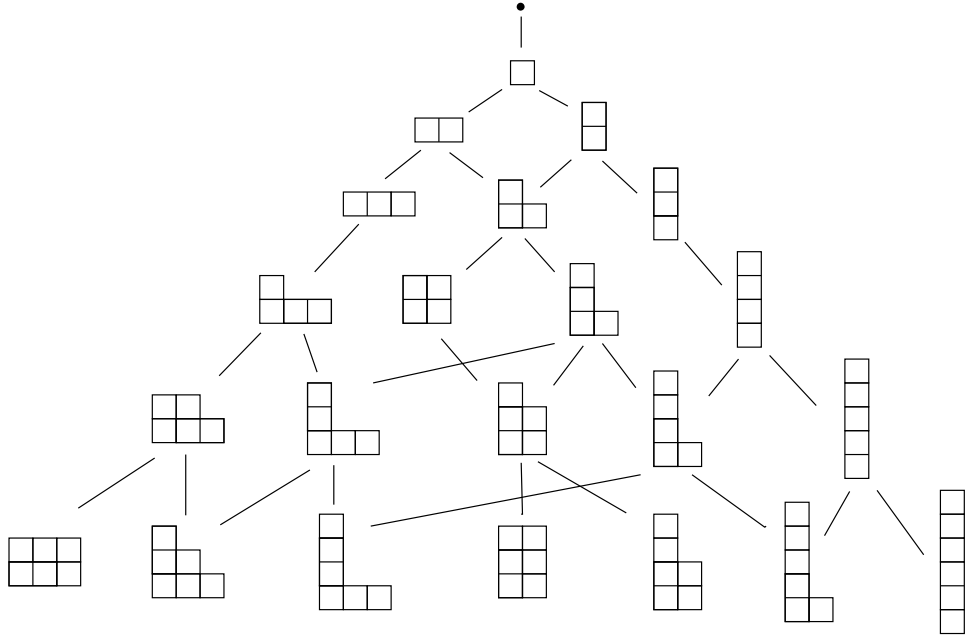}
\end{center}
\caption{$3$-analogue of the Young Lattice}
\label{klattice}
\end{figure}

When $\ell=1$ in \ref{pieri}, we are led to a $k$-analogue 
of the Young lattice.  This is the poset whose elements 
are $k$-bounded partitions and whose Hasse diagram is obtained by linking
an element $\lambda$ to every $\mu\in E_{\lambda,1}^{(k)}$.  
In Figure~\ref{klattice}, we illustrate the poset obtained for 
degree 6 with $k=3$.  Moreover, the number of paths in this poset 
joining the empty partition to the partition $\lambda$ is 
simply the number of summands in \ref{kostkafoulkes1} 
when $\mu=1^{|\lambda|}$, namely $K_{\lambda,1^{|\lambda|}}^{(k)}(1)$.
An analogous observation can be made for a general $\mu$.

Central to our research is the observation that not all of the atoms 
need to be constructed using \ref{defatom1}.  In fact, for each $k$ 
there is a distinguished ``irreducible'' subset of atoms of $V_k$ from 
which all successive atoms may be constructed simply by 
applying certain generalized promotion operators.
To be more precise, let a partition $\mu$ with no more than 
$i$ parts equal to $k-i$ be called $k$-irreducible
(note, there are $k!$ such partitions).
Thus, any $k$-bounded partition can be obtained by the partition
rearrangement of the parts of a $k$-irreducible partition 
and a sequence of $k$-rectangles, partitions of
the form $(\ell^{k+1-\ell})$ for $\ell=1,\dots,k$.
This given, we let the collection of $k$-irreducible atoms 
be only the atoms indexed by $k$-irreducible partitions.
This suggests that there are certain generalized
promotion operators indexed by $k$-rectangles
yielding that every atom may be written in the form
\begin{equation}
A^{(k)}_\lambda[X;t]
=
t^c
\digamma\left(
\mathbb B_{R_1}
\mathbb B_{R_2}
\ldots
\mathbb B_{R_\ell}
\mathbb A^{(k)}_\mu
\right) \, ,
\label{atirr}
\end{equation}
where $\mu$ is a $k$-irreducible partition,
$R_1,\ldots,R_\ell$ are certain $k$-rectangles,
and $c\in\mathbb N$.

Again we find it interesting to consider the
case $t=1$.  First, since the Hall-Littlewood polynomials 
at $t=1$ are simply
\begin{equation}
H_{(\mu_1,\ldots,\mu_n)}[X;1]
=
h_{\mu_1}[X]\cdots h_{\mu_n}[X]
\,,
\end{equation}
we see that $V_k$ reduces to the polynomial
ring $V_k(1)=\mathbb Q[h_1,\ldots,h_k]$.
Further, since the construction in \ref{atirr} is 
simply multiplication by Schur functions when $t=1$,
\begin{equation}
A^{(k)}_\lambda[X;1]
=
S_{R_1}[X]
S_{R_2}[X]
\ldots
S_{R_\ell}[X]
A^{(k)}_\mu[X;1]
\,,
\label{surmul}
\end{equation}
it thus follows that $k$-irreducible atoms constitute a natural basis
for the quotient $V_k(1)/{\cal I}_k$, where
${\cal I}_k$ is the ideal generated by 
Schur functions indexed by $k$-rectangles.
In fact, the irreducible atom basis offers
a very beautiful way to carry out 
operations in this quotient ring:  first work in
$V_k(1)$ using atoms and then
replace by zero all atoms indexed by partitions which are
not $k$-irreducible.

We shall examine our $k$-analogue of the Young lattice 
restricted to $k$-irreducible partitions.  
Figure~\ref{figindecomp} gives the case $k\!=\!3$ and $k\!=\!4$, where 
vertices denote irreducible atoms rather than partitions.  
Since it can be shown that 
the collection of monomials of the form
$\{ h_1^{\epsilon_1}, h_2^{\epsilon_2}, \ldots , h_k^{\epsilon_k} \}_{
0\leq \epsilon_i\leq k-i}$
provides a basis for the quotient $V_k(1)/{\cal I}_k$,
it follows that the Hilbert series
$F_{V_k(1)/{\cal I}_k}(q)$
of this quotient, as well as the rank generating function of the 
corresponding poset, is given by 
\begin{equation}
F_{V_k(1)/{\cal I}_k}(q)
=
\prod_{i=1}^{k-1} 
\left(
1+q^i+q^{2i}+\cdots + q^{(k-i)i}
\right) \, .
\label{HS}
\end{equation}

Finally, we shall make connections between our work and contemporary 
research in this area.  We discovered that tableaux manipulations identical
to ours have been used for a different purpose in \cite{[S2],[S3],[S1]}.
In particular, 
certain cases of the generalized Kostka polynomials 
can be expressed in our notation as 
\begin{equation}
\digamma \left( \mathbb B_{R_1} \cdots \mathbb B_{R_\ell} \, \mathbb H_0 \right) \,, 
\label{prodrec1}
\end{equation}
where $R_1,\ldots,R_\ell$ is a sequence of rectangles whose
concatenation is a partition \cite{[S3]}.  When $R_1,\ldots,R_\ell$ is a sequence of 
$k$-rectangles this is simply the case $\mu=\emptyset$ in \ref{atirr}. 
Thus, it is again apparent that an integral part of our work lies in 
the $k$-irreducible atoms, without which the atoms in general could not 
be constructed.

Furthermore, it is known that these generalized Kostka polynomials can be
built from the symmetric function operators $B_R$ introduced in \cite{[ZS]}.
The connection we have made with our atoms and these polynomials
thus suggest that any atom can be obtained by applying a succession of operators
$B_R$ indexed by $k$-rectangles to a given irreducible atom $A_{\mu}^{(k)}[X;t]$:
\begin{equation}
A_{\lambda}^{(k)}[X;t] = 
t^c B_{R_1} B_{R_2} \cdots  B_{R_\ell} 
A_{\mu}^{(k)}[X;t] \, ,
\quad\text{where}\; c\in\mathbb N
\,.
\end{equation}
Note this is a symmetric function analogue of \ref{atirr} 
and specializes to \ref{surmul} when $t=1$.

\begin{acknow}
We give our deepest thanks to Adriano Garsia for all his
time and effort helping us articulate our ideas.  L. Lapointe would
also 
like to thank Luc Vinet for his support and the helpful discussions.
Our research depended on the use of ACE \cite{[V]}.
\end{acknow}

\interligne{0.4}

\tableofcontents

\interligne{1.0}

\section{Background}

\subsection{Symmetric function theory}

Here, symmetric functions are indexed by partitions, or sequences of
non-negative integers $\lambda =(\lambda_1,\lambda_2,\ldots)$ with 
$\lambda_1 \ge \lambda_2 \ge \dots$. 
The order of $\lambda$ is $|\lambda| = \lambda_1 + \lambda_2 + \dots$, 
the number of non-zero parts in $\lambda$ is denoted $\ell(\lambda)$,
and  $n(\lambda)=\sum_{i} (i-1) \lambda_i$.
We use the dominance order on partitions with $|\lambda|=|\mu|$, where 
$\lambda\leq\mu$ when $\lambda_1+\cdots+\lambda_i\leq
\mu_1+\cdots+\mu_i$ for all $i$.  For two partitions $\lambda$
and $\mu$, $\lambda \cup \mu$ denotes the partition rearrangement of
the parts of $\lambda$ and $\mu$.

Every partition $\lambda$ may be associated to 
a Ferrers diagram with $\lambda_i$ lattice squares in the $i^{th}$
row, from the bottom to top.  For example,
\begin{equation}
\lambda\,=\,(4,3,1)\,=\,
{\tiny{\tableau*[scY]{ \cr & & \cr & & & }}} \, .
\end{equation}
For each cell $s=(i,j)$ in the diagram of $\lambda$, let
$\ell'(s), \ell(s), a(s)$, and $a'(s)$ be respectively the number of
cells in the diagram of $\lambda$ to the south, north, east, and west
of the cell $s$.
The transposition of a diagram associated to $\lambda$ 
with respect to the main diagonal gives the conjugate partition
$\lambda'$.  For example, the conjugate of (4,3,1) is
\begin{equation}
\lambda' \,=\,
{\tiny{\tableau*[scY]{ \cr & \cr & \cr & &}}} 
\,=\,(3,2,2,1)\,.  
\end{equation}
\smallskip

A skew diagram $\mu/\lambda$, for any partition $\mu$
containing the partition $\lambda$, is the diagram obtained by 
deleting the cells of $\lambda$ from $\mu$.
The thick frames below represent (5,3,2,1)/(4,2).
\begin{equation}
{\tiny{\tableau*[scY]{ \tf \cr \tf & \tf \cr  & & 
\tf \cr & & & & \tf }}} \, .
\end{equation}

We employ the notation of $\lambda$-rings, needing only the formal ring of symmetric
functions $\Lambda$ to act on the ring of rational functions in $x_1,\dots,x_N,q,t$,
with coefficients in $\mathbb Q$.  
The action of a power sum $p_i$ on a rational function is, by definition,
\begin{equation}
p_{i} \left[ \frac{\sum_{\alpha} c_{\alpha} u_{\alpha} 
 }{ \sum_{\beta} d_{\beta} v_{\beta} } \right]
 =\frac{\sum_{\alpha} c_{\alpha} u_{\alpha}^i }{ \sum_{\beta} d_{\beta} v_{\beta}^i}, 
\label{pleth}
\end{equation}
with $c_{\alpha},d_{\beta} \in \mathbb Q$ and $u_{\alpha}, v_{\beta}$ 
monomials in $x_1,\dots,x_N,q,t$.  
Since the ring $\Lambda$ is generated by power sums, $p_i$,
any symmetric function has a unique expression in terms of $p_i$, 
and \ref{pleth} extends to an action of $\Lambda$ on rational functions.  
In particular $f[X]$, 
the action of a symmetric function 
$f$ on the polynomial $X=x_1+\cdots+x_N$, 
is simply $f(x_1,\ldots,x_N)$.  
\smallskip

We recall that the Macdonald scalar product,  
$\langle \ , \ \rangle_{q,t}$, 
on $\Lambda \otimes \mathbb Q(q,t)$ is defined by setting 
\begin{equation}
\langle p_\lambda[X], p_\mu[X] \rangle_{q,t}
        =\delta_{\lambda \mu } \, z_\lambda \prod_{i=1}^{\ell(\lambda)} \frac{1-q^{\lambda_i}}{1-t^{\lambda_i}} = \delta_{\lambda \mu } \,  z_\lambda \, p_{\lambda}\left[\frac{1-q}{1-t}\right] \,,
\end{equation}
where for a partition $\lambda$ with $m_i(\lambda)$ parts equal to $i$,
we associate the number
\begin{equation} \label{1}
z_\lambda
        = 1^{m_1} m_1 !  \, 2^{m_2} m_2! \dotsm
\end{equation}
The Macdonald integral forms  $J_\lambda [X; q,t]$ are then uniquely characterized
\cite{[Ma]} by
\begin{align}
  \mathrm{(i)} \ &  \langle J_\lambda, J_\mu \rangle_{q,t} = 0, \qquad \text{if } \lambda \ne \mu , \\
  \mathrm{(ii)} \ &  J_\lambda[X;q,t] = \sum_{\mu \le \lambda} v_{\lambda\mu }(q,t) 
S_\mu[X] , \\
  \mathrm{(iii)} \ & v_{\lambda\lambda}(q,t)= 
\prod_{s \in \lambda} (1-q^{a(s)} t^{\ell (s)+1}),
\end{align}
where $S_{\mu}[X]$ is the usual 
Schur function and $v_{\lambda \mu}(q,t) \in \mathbb Q(q,t)$.

\subsection{Tableaux combinatorics} \label{subsectab}

$\mathcal A^{*}$ denotes the free monoid generated by 
the alphabet $\mathcal A = \{1,2, \dots \}$
and $\mathbb Q [\mathcal A^{*}]$ is the free algebra of $\mathcal A$.  
Elements of $\mathcal A^{*}$ are called words
and for $E$ a subset of $\mathcal A$, $w_E$ denotes the subword 
obtained by removing from $w$ all the letters not in $E$.
The degree of a word $w$ is denoted $|w|$ 
and if $w$ has $\rho_1$ ones, $\rho_2$ twos, $\ldots$, and
$\rho_m$ $m$'s, then the evaluation of $w$ is 
$(\rho_1,\ldots,\rho_m)$.   
For example, $w = 1 3 1  3 3 2$ has degree 6 and evaluation (2,1,3). 
A word $w$ of degree $n$ is standard iff its evaluation is $(1,\dots,1)$.
Recall that a word $w$ is Yamanouchi in the letters 
$a_1\!<\!\dots\!<\!a_{h}$ if it is such that for every factoring $w=u v$, 
$v$ contains more $a_i$ than $a_j$ for all $i\!<\!j$.  

The plactic monoid on the alphabet $\mathcal A$ is the quotient 
$\mathcal A^{*}/\equiv$, where $\equiv$ is the congruence 
generated by the Knuth relations \cite{[Kn]}
defined on three letters $a,b,c$ by
\begin{equation}
\begin{split}
a \, c \, b  &\equiv c \, a \, b \qquad ( a \leq b < c) \, , \\ 
b \, a \, c   &\equiv b \, c \, a  \qquad ( a <  b \leq c) \, . \\ 
\end{split} 
\label{Knuth}
\end{equation}
Two words $w$ and $w'$ are said to be Knuth equivalent iff $w \equiv w'$.

In this paper, a tableau 
is a filling of a Ferrers diagram with positive integer 
entries that are nondecreasing in rows
and increasing in columns: 
\begin{equation}
\sstab = {\tiny{\tableau*[scY]{
6 & 7 \cr
4 & 4 & 5 \cr
1 & 1 &1 &2 & 3 }}} \, .
\label{tab}
\end{equation}
The word $w$ obtained by reading the entries of a tableau 
from left to right and top to bottom is said to be a tableau word, 
or simply a tableau.  Our example shows that
$w=6744511123$ is a tableau
with evaluation (3,1,1,2,1,1,1).  
A standard tableau $\stab$ is a 
tableau of evaluation $(1,1,\dots,1)$.  For example,\begin{equation}
\stab = {\tiny{\tableau*[scY]{
7  \cr
 4 & 6 \cr
1 & 2 &3 &5 }}} \, .
\label{stab}
\end{equation}

The transpose of a standard tableaux $\stab^t$ is defined in the same 
manner as the transpose of a Ferrers diagram.
With $\stab$ as given in \ref{stab}, we have
\begin{equation}
\stab^t = {\tiny{\tableau*[scY]{ 5 \cr 3 \cr 2 & 6 \cr 1 & 4 & 7  }}} \, .
\end{equation}
Since the transpose of a tableau is assured to be
a tableau only when the original tableau is standard,
this definition is valid only for standard tableaux.

We assume readers are familiar with the Robinson-Schensted correspondence 
\cite{[Rob],[Sch]},
\begin{equation}
w \longleftrightarrow \bigl(P(w),Q(w) \bigr) \, ,
\label{RS}
\end{equation}
providing a bijection between a word $w$ and a pair 
of tableaux $\bigl(P(w),Q(w)\bigr)$, where $P(w)$ is the
only tableau Knuth equivalent to $w$ and $Q(w)$ is a standard tableau.

The ring of symmetric functions is embedded into the plactic
algebra by sending the Schur function $S_\lambda$ to the sum of all tableaux
of shape $\lambda$ \cite{[Tab],[Lo]}.  The commutativity of the product 
$S_\lambda S_\mu\! \equiv\! S_\mu S_\lambda$ 
thus implies bijections among tableaux.
In particular, we can define the following action of
the symmetric group on words \cite{[LS1]}.  
The elementary transposition $\sigma_i$ permutes degrees in 
$i$ and $i+1$. Given a word $w$ of evaluation
$(\rho_1,\dots,\rho_m)$, let $u$ denote
the subword in letters $a=i$ and $b=i+1$.
The action of the transposition $\sigma_i$ affects only the subword
$u$ and is defined as follows:
pair every factor $b\, a$ of $u$, and let $u_1$ be the subword of $u$ 
made out of the unpaired letters.  Pair every factor $b\,a$ of $u_1$, and 
let $u_2$ be the subword made out of the unpaired letters.  Continue
in this fashion as long as possible.
When all factors $b\,a$ are paired and unpaired letters of $u$
are of the form $a^r b^s$, $\sigma_i$ sends $a^rb^s\to a^sb^r$.
For instance, to obtain the action of $\sigma_2$ on $w=123343222423$, 
we have $u=w_{\{2,3\}}=233322223$, and the pairings are
\begin{equation}
2\Bigl(3\bigl(3(32)2\bigl)2 \Bigl)23 \, ,
\end{equation}
which means that
$\sigma_2 u=233322233$ and $\sigma_2 w=123343222433$.  
It is verified in \cite{[LS1]} that the $\sigma_i$'s obey the Coxeter relations
and thus provide an action of the symmetric group on words.

We use the  notion of charge \cite{[Lo],[LS2]} defined by writing a word $w$,
with evaluation given by a partition,
counterclockwise on a circle with a * separating 
the end of the word from its beginning and then 
summing the labels that are obtained by the following procedure: 
Let $\ell=0$.
Moving clockwise from *, we label with $\ell$,
the first occurrence of letter $1$, then the first occurrence of 
letter $2$ following this $1$, then the first occurrence of letter 
$3$ following this $2$, etc, with the condition that
each time the * is passed, the label $\ell$ is increased to $\ell+1$.
Once each of the letters $1,2,3,\ldots$ have been labeled,
we repeat this procedure on the unlabeled letters, again
starting at the * with $\ell=0$. The process ends when all letters have 
been labeled.

We can define charge on a word
$w$ whose evaluation is not a partition by first permuting the
evaluation to a partition using $\sigma$,
and then taking the charge of $\sigma w$.
\begin{figure}[htb]
\begin{center}
\epsfig{file=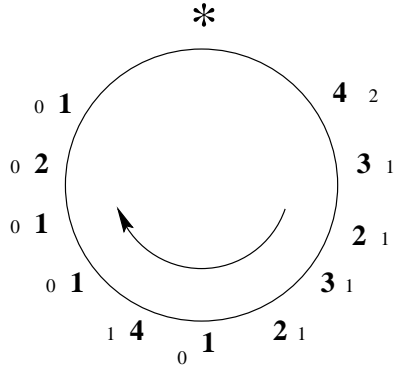}
\end{center}
\caption{Charge of $w=12114123234$}
\label{charge}
\end{figure}
Figure~\ref{charge} shows that
$\charge(12114123234)=0+0+0+0+1+0+1+1+1+1+2=7$.

The definition of charge gives that
a tableau of shape and evaluation $\mu$
has charge 0 and thus the combinatorial construction 
for the Hall-Littlewood polynomials
(\ref{halldef}) implies 
\begin{equation}
H_{\mu}[X;t] = S_{\mu}[X]+ \sum_{\lambda > \mu}
K_{\lambda \mu}(t) S_{\lambda}[X] \, .
\label{HLF}
\end{equation}

\section{Definition of $A_{\lambda}^{(k)}[X;t]$}
\label{secatom}
Our main contribution is the method for constructing 
new families of functions whose significance has been 
outlined in the introduction.
The characterization is similar to 
the combinatorial definition of Hall-Littlewood 
polynomials using the set of tableaux that arises
from a recursive application of promotion operators.
Our families also correspond to a set of tableaux 
generated by the promotion operators, but here we introduce 
new operators $\mathbb P_{\lambda^{\to k}}$ to eliminate
undesirable elements.  To be precise, we now define the 
operators involved in our construction.

The promotion operators are defined
on a tableau $T$ with evaluation 
$(\lambda_1,\ldots,\lambda_m)$ by
\begin{equation}
\mathbb B_{r}( 
T) =  \sigma_{1} \cdots \sigma_{m} \,  
\mathbb R_{r} \, T 
\, ,
\end{equation}
where $\mathbb R_r$ adds a horizontal $r$-strip 
of the letter $m+1$ to $T$ in all possible ways.  For example,
\begin{equation*}
\mathbb B_{3} \, {\tiny{\tableau*[scY]{ 2 & 2\cr 1 & 1}}} 
\, =\, \sigma_{1} \, \sigma_2 \,
\mathbb R_{3} \, {\tiny{\tableau*[scY]{ 2 & 2\cr 1 & 1}}}
\, = \,   
\sigma_{1} \sigma_2 \left(\,  
{\tiny{\tableau*[scY]{ 2 & 2 \cr 1 & 1 & 3 & 3 & 3}}}
+{\tiny{\tableau*[scY]{3\cr 2 & 2 \cr 1 & 1 & 3 & 3}}}
+{\tiny{\tableau*[scY]{ 3 & 3\cr 2 & 2 \cr 1 &1 & 3}}}
\,
\right)
\, = \,
{\tiny{\tableau*[scY]{ 2 & 2 \cr 1 & 1 & 1 & 3 & 3}}}
+{\tiny{\tableau*[scY]{3\cr 2 & 2 \cr 1 & 1 & 1 & 3}}}
+{\tiny{\tableau*[scY]{ 3 & 3\cr 2 & 2 \cr 1 &1 & 1}}}
\,.
\end{equation*}
Note that the action of $\sigma$ implies that
the resulting tableaux have evaluation
$(r,\lambda_1,\ldots,\lambda_m)$.
While our construction relies on these operators,
they generate certain unwanted tableaux.
We now present the concepts needed to obtain
operators that filter out such elements.  

The {\it main hook-length} of a partition $\lambda$, 
$h_M(\lambda)$, is the hook-length of the cell $s=(1,1)$
in the diagram associated to $\lambda$.  That is
\begin{equation}
h_M(\lambda) = \ell(s)+a(s)+1 = \lambda_1+\lambda_1'-1 = \lambda_1 + \ell(\lambda) -1 \, .
\end{equation}
For example, if $\lambda=(4,3,1)$, then 
\begin{equation}
h_M(\lambda)=
{\tiny{\tableau*[scY]{ \bullet \cr \bullet & & \cr \bullet & \bullet & \bullet & \bullet}}}
= 6 \;  .
\end{equation}
Any $k$-bounded partition $\lambda$ 
can be associated to a sequence of partitions 
called the $k$-{\it split},
$\lambda^{\to k}= (\lambda^{(1)},\lambda^{(2)},\dots,\lambda^{(r)})$.  
The $k$-{split} of $\lambda$ is obtained by 
dividing $\lambda$ (without changing the order of its entries) into partitions $\lambda^{(i)}$ 
where $h_M(\lambda^{(i)})= k$, $\forall \, i \neq r$.  For example, 
$(3,2,2,2,1,1)^{\to 3}=\bigl((3),(2,2),(2,1),(1)\bigr)$ .
Equivalently, we horizontally cut the diagram of $\lambda$
into partitions $\lambda^{(i)}$ where $h_M(\lambda^{(i)})= k$. In our
example this gives
\begin{equation}
{\tiny{\tableau*[scY]{ \cr \cr & \cr & \cr & \cr & & \cr}}} \quad 
\begin{matrix}
& {\tiny{\tableau*[sbY]{  \cr}}} \hfill \\
& {\tiny{\tableau*[sbY]{  \cr & \cr}}} \hfill \\ 
\longrightarrow & {\tiny{\tableau*[sbY]{ & \cr  & \cr}}} \hfill \\ 
& {\tiny{\tableau*[sbY]{  & & \cr}}} \hfill 
\end{matrix} \, .
\end{equation}
Note, the last partition in the sequence 
$\lambda^{\to k}$ may have main hook-length less than $k$.
As $k$ increases, the $k$-{split} of $\lambda$ will contain 
fewer partitions.  For $k=4$,
\begin{equation}
{\tiny{\tableau*[scY]{ \cr \cr & \cr & \cr & \cr & & \cr}}} \quad 
\begin{matrix}
& {\tiny{\tableau*[sbY]{  \cr}}} \hfill \\ 
\longrightarrow & {\tiny{\tableau*[sbY]{\cr  & \cr  & \cr}}} \hfill \\ 
& {\tiny{\tableau*[sbY]{ & \cr & & \cr}}} \hfill 
\end{matrix} \qquad {\rm{or}} \qquad
 (3,2,2,2,1,1)^{\to 4}=\bigl((3,2),(2,2,1),(1)\bigr) \, . 
\end{equation}
When $k$ is big enough 
($h_M(\lambda)\leq k$), then
$\lambda^{\to k}\!=\!(\lambda)$.  
i.e., $(3,2,2,2,1,1)^{\to 8}\!=\!\bigl((3,2,2,2,1,1) \bigr)$.

Let $T$ be a given tableau whose shape contains $\lambda$.  
We shall denote by $T_{\lambda}$, the subtableau of $T$ of shape $\lambda$.
Let $U$ be the skew tableau obtained by removing
$T_{\lambda}$ from $T$, let $T_1$ be the tableau contained in 
the first $\ell(\lambda)$ rows of 
$U$, and $T_2$ be the portion of $U$ that is above the 
$\ell(\lambda)$ rows.
Let us denote by $T_1 T_2$ the skew tableau obtained by juxtaposing
$T_1$ to the northwest corner of $T_2$, and by $\overline{T}$ the unique
tableau which is Knuth equivalent to $T_1T_2$.
For instance, in the figure below $\lambda=(3,2,1,1)$,
the skew tableau with empty cells is $T_1$, the 
tableau with bullets is $T_2$, the middle diagram
is $T_1T_2$, and the right diagram is a possible shape for $\overline{T}$.
\begin{equation}
{\tiny{\tableau*[scY]{ \bullet \cr \bullet & \bullet \cr \bullet &  \bullet & \bullet & 
\bullet \cr \fl &   & & \cr \fl &  &   & & \cr \fl&\fl  &  & & & \cr \fl & \fl & \fl & & & & & }}}
\; \longrightarrow 
\;
{\tiny{\tableau*[scY]{ & & \cr  & & & \cr \bl & & & & \cr \bl & \bl &  & & & & \cr 
\bl & \bl & \bl & \bl & \bl & \bl & \bl & \bullet  \cr
\bl & \bl & \bl & \bl & \bl & \bl & \bl & \bullet & \bullet \cr
\bl & \bl & \bl & \bl & \bl & \bl & \bl & \bullet & \bullet & 
\bullet & \bullet \cr}}} \; \equiv \; 
{\tiny{\tableau*[scY]{ \cr \cr & \cr & & \cr & & \cr & & & & & \cr & & & & & &}}} 
\, . 
\end{equation}
This construction permits us to define an operation on tableau,
$\mathbb K_\lambda$, called {\it $\lambda$-katabolism}.
\begin{equation*}
\mathbb K_\lambda (
T)=
\begin{cases}
\overline{T} & \text{ if } \lambda\subseteq \text{shape}(T)\\
0 & \text{ otherwise } 
\end{cases}
\,.
\end{equation*}
For example, the $(2,1)${-katabolism} of $T=9472581236$ is
\begin{equation} 
\mathbb K_{(2,1)}\,\,
{{\tableau*[scY]{ 9 \cr   4&7\cr \tf 2& 5 &8 \cr \tf 1& \tf 2&  3&6}}}  
\quad \longrightarrow \quad
{{\tableau*[scY]{  5 & 8 \cr \bl &  3 & 6 \cr \bl & \bl & \bl &  9 \cr  \bl & \bl & \bl & 4 & 7 \cr }}} \quad \equiv \quad
{{\tableau*[scY]{ 8\cr 5 & 6 & 9  \cr 3 & 4 & 7}}} 
\, \, \, .
\label{katex}
\end{equation}
Note that $\lambda$-katabolism was also introduced in \cite{[S1]}
and for the case that $\lambda$ is a row, in \cite{[LS1]}.

Let $S(\lambda)$ be the set of $\lambda$-shaped tableaux with evaluations 
$(0^m,\lambda_1,\lambda_2,\ldots)$, for $m\in\mathbb N$.
For $\lambda=(3,2,2)$, we have
\begin{equation}
S\left(\,{\tiny{\tableau*[scY]{&  \cr &\cr && }}}\, \right) = \left \{ \, 
{\tiny{\tableau*[scY]{3 & 3 \cr 2&2\cr 1&1&1 }}} \, \, , \, \, 
{\tiny{\tableau*[scY]{4 & 4 \cr 3&3\cr 2&2&2 }}} \, \, , \, \, 
{\tiny{\tableau*[scY]{5 & 5 \cr 4&4\cr 3&3&3 }}} \, \, , \, \, 
{\tiny{\tableau*[scY]{6 & 6 \cr 5&5\cr 4&4&4 }}} \, \, , \, \, 
\dots \, \right \} \, \, .
\end{equation}
This given, the restricted $\lambda${-katabolism} 
$\mathbb {\overline K}_\lambda$ is
defined by setting
\begin{equation}
\mathbb {\overline K}_\lambda ( 
T)= 
\begin{cases}
\mathbb K_\lambda(T) 
& {\text{if }\,  T_{\lambda} \in S(\lambda) } \\
0 & {\text{otherwise}} \\
\end{cases} \, . 
\end{equation}
For example, $\mathbb {\overline K}_{(2,1)}$ on the tableau in \ref{katex}
is zero, whereas
\begin{equation}
\mathbb {\overline K}_{(2,1)}\,\;
{{\tableau*[scY]{ 9 \cr   4&7\cr \tf 2& 5 &8 \cr \tf 1& \tf 1&  3&6}}}  
\quad = \quad
{{\tableau*[scY]{  5 & 8 \cr \bl &  3 & 6 \cr \bl & \bl & \bl &  9 \cr  \bl & \bl & \bl & 4 & 7 \cr }}} \quad \equiv \quad
{{\tableau*[scY]{ 8\cr 5 & 6 & 9  \cr 3 & 4 & 7}}} \, \, \, .
\end{equation}

For a sequence of partitions 
$S= (\lambda^{(1)},\lambda^{(2)},\ldots,\lambda^{(\ell)})$, we define
the filtering operator $\mathbb P_{S}$ using 
the succession of restricted katabolisms 
$\overline{\mathbb K}_{\lambda^{(\ell)}}\cdots
\overline{\mathbb  K}_{\lambda^{(1)}}$, 
\begin{equation}
\mathbb P_{S} ( T )=
\begin{cases}
T & {\text{if $\overline{\mathbb K}_{\lambda^{(\ell)}}
\cdots \overline{\mathbb K}_{\lambda^{(1)}}(T)= 
\mathbb H_{0} $ }}\\
0 & {\text{otherwise}} \\
\end{cases} \, \, .
\label{FilOp}
\end{equation}
In fact, we only consider the case where $S$ is 
the sequence of partitions given by $\lambda^{\to k}$.
\begin{property} 
The filtering operators $\mathbb P_{\lambda^{\to k}}$ satisfy the following properties :

a. For $T\in S(\lambda)$, we have 
$\mathbb P_{\lambda^{\to k}} T = T$ 
for all $k$ such that $\lambda$ is bounded by $k$.

b. For $U$ a tableau of $|\lambda|$ letters such that
$U\not\in S(\lambda)$,
$\mathbb P_{\lambda^{\to k}} U = 0$
for all $k\geq h_M(\lambda)$. 
\label{FilProp}
\end{property}
\noindent {\it Proof.} \quad 
By definition \ref{FilOp}, 
we must show
$\overline{\mathbb K}_{\lambda^{(\ell)}}\cdots
\overline{\mathbb K}_{\lambda^{(1)}}T=\mathbb H_0$,
for $\lambda^{\to k}=( \lambda^{(1)}, \ldots, \lambda^{(\ell)})$.
Recall that $\overline{\mathbb K}_{\lambda^{(1)}} T$ 
acts by extracting the bottom $\ell(\lambda^{(1)})$ rows of $T$
and inserting into the remainder, 
any entries not in $T_{\lambda^{(1)}}\in S(\lambda^{(1)})$.
Since the bottom $\ell(\lambda^{(1)})$ rows of $T\in S(\lambda)$
are exactly $T_{\lambda^{(1)}}\in S(\lambda^{(1)})$,
the katabolism $\overline{\mathbb K}_{\lambda^{(1)}} T$ 
simply removes the bottom $\ell(\lambda^{(1)})$ rows of $T$.
By iteration, we obtain the empty tableau.

For $(b)$, the condition that $k$ is large implies that 
$\lambda^{\to k}=(\lambda)$.  
It thus suffices to show that 
$\overline{\mathbb K}_{\lambda}(U)=0$.
Now $\overline{\mathbb K}_{\lambda}$ acts first by extracting 
from $U$, the subtableau $U_\lambda\in S(\lambda)$.
If $U$ is of shape $\lambda$ then
$U_\lambda=U\not\in S(\lambda)$.
If $U$ is not of shape $\lambda$,
since $U$ has degree $|\lambda|$,
then $U_\lambda$ does not exist.
Therefore we have our claim.
\hfill 
$\endprf$

These filtering operators are those required in the characterization of our 
families of functions.
We thus have the tools to recursively define the 
central object in our work, the {\it super atom} of 
shape $\lambda$ and level $k$, $\mathbb A_{\lambda}^{(k)}$. 
\begin{definition}
Let $\mathbb A_{0}^{(k)}$ be the empty tableau.
The super atom of a
$k$-bounded partition $\lambda$ is
\begin{equation}
\mathbb A_{\lambda}^{(k)} = \mathbb P_{\lambda^{\to k}} \, 
\mathbb B_{\lambda_1} \, \left( \mathbb A_{(\lambda_2,\lambda_3,\dots)}^{(k)} \right) \, .
\end{equation}
\label{defatom}
\end{definition}
\noindent 
For example, given that we know the super atom
\begin{equation}
A_{1,1,1,1}^{(3)} \,= \, {\tiny{\tableau*[scY]{ 4 \cr 3 \cr 2 \cr 1 }}} +
{\tiny{\tableau*[scY]{  3 \cr 2 \cr 1 & 4 }}} 
\label{At1111}
\end{equation}
we can obtain $\mathbb A_{2,1,1,1,1}^{(3)}$
by first acting with the rectangular operator $\mathbb B_{2}$
on \ref{At1111},
\begin{equation}
\mathbb B_{2} \, 
\left( \mathbb A_{1,1,1,1}^{(3)} \right) \,=\, \mathbb B_{2} \, 
\left( \, {\tiny{\tableau*[scY]{ 4 \cr 3 \cr 2 \cr 1 }}} +
{\tiny{\tableau*[scY]{  3 \cr 2 \cr 1 & 4 }}} \,   \right) \, = \,  
{\tiny{\tableau*[scY]{ 3 \cr 2 \cr 1 & 1 & 4 & 5 }}} + 
{\tiny{\tableau*[scY]{  3 \cr 2 & 5 \cr 1 & 1 & 4 }}} + 
{\tiny{\tableau*[scY]{ 5 \cr 3 \cr 2 \cr 1 & 1 & 4 }}} +
{\tiny{\tableau*[scY]{ 4 \cr 3 \cr 2 & 5 \cr 1 & 1 }}} + 
{\tiny{\tableau*[scY]{ 4 \cr 3 \cr 2 \cr 1 & 1 & 5}}} +
{\tiny{\tableau*[scY]{ 5 \cr 4 \cr 3\cr 2 \cr 1& 1 }}} \, ,
\label{exatom}
\end{equation}
and then to these tableaux, applying the operator $\mathbb P_{\lambda^{\to 3}}$,
where $\lambda^{\to 3} = \bigl( (2,1),(1,1,1) \bigr)$;
\begin{equation}
\mathbb A_{2,1,1,1,1}^{(3)} \, = \, 
\mathbb P_{((2,1),(1,1,1))} \, \left( 
\mathbb B_{2} \, 
 \mathbb A_{1,1,1,1}^{(3)} \right)\,=\, 
{\tiny{\tableau*[scY]{  3 \cr 2 & 5 \cr 1 & 1 & 4 }}} + 
{\tiny{\tableau*[scY]{ 4 \cr 3 \cr 2 & 5 \cr 1 & 1 }}} + 
{\tiny{\tableau*[scY]{ 4 \cr 3 \cr 2 \cr 1 & 1 & 5}}} +
{\tiny{\tableau*[scY]{ 5 \cr 4 \cr 3\cr 2 \cr 1& 1 }}} \, .
\end{equation}

Our method for constructing the super atoms allows the 
derivation of several natural properties.
In particular, these properties generally arise as
the consequence of those held by the promotion 
and filtering operators. 
\begin{property}
For all $k$-bounded partitions $\lambda$,
we have
\begin{equation}
\mathbb A_\lambda^{(k)}\subseteq \mathbb H_\lambda \, .
\end{equation}
\label{atominhall}
\end{property}
\noindent {\it Proof.} \quad 
For $\lambda=(\lambda_1,\dots,\lambda_m)$, recall from \ref{hallrec}  
that
\begin{equation}
\mathbb H_{\lambda} = 
\mathbb B_{\lambda_1} \cdots \mathbb B_{\lambda_{m}} 
\,  \mathbb H_{0} \, .
\end{equation}
On the other hand, 
following from Definition~\ref{defatom}, we have 
\begin{equation}
\mathbb A_{\lambda}^{(k)} = 
\mathbb P_{\lambda^{\to k}} \, \mathbb B_{\lambda_1} \cdots 
\mathbb P_{(\lambda_m)^{\to k}}  \,  \mathbb B_{\lambda_{m}} 
\,  \mathbb A_{0}^{(k)} \, .
\end{equation}
Since $\mathbb A_{\lambda}^{(k)}$ 
is distinguished from $\mathbb H_{\lambda}$ 
only by acting with a filtering operator
after each application of a $\mathbb B_{\ell}$ operator, 
we have that every tableau in 
$\mathbb A_{\lambda}^{(k)}$ is also in $\mathbb H_{\lambda}$.
\hfill 
$\endprf$

\begin{property}
Let $T$ be the tableau of shape and evaluation $\lambda$.
The super atoms satisfy:

i.  $T\in\mathbb A_\lambda^{(k)}$  for any $k$ such that $\lambda$ is bounded by $k$.

ii.  $\mathbb A_\lambda^{(k)} = T$  for $ k\geq h_M(\lambda)$. 
\label{SAProp}
\end{property}
\noindent {\it Proof.} \quad 
($i$):  Recall that
$\mathbb A_{\lambda}^{(k)}
= \mathbb P_{\lambda^{\to k}} \mathbb B_{\lambda_1}
\mathbb A_{\lambda_2,\ldots,\lambda_\ell}^{(k)}$.
Assume by induction that 
$U\!\in\mathbb A_{\lambda_2,\ldots,\lambda_\ell}^{(k)}$
where $U$ has shape and evaluation 
$(\lambda_2,\ldots,\lambda_\ell)$.
$\mathbb R_{\lambda_1}U$ produces a sum of tableaux,
one being the tableau of shape $(\lambda_1,\lambda_2,\ldots,\lambda_\ell)$
which is then sent to the tableau $T$ of shape and evaluation $\lambda$
under the action of the symmetric group.
It thus suffices to show that $T$ is not eliminated 
by $\mathbb P_{\lambda^{\to k}}$ for all $k$.
This is shown in Property~\ref{FilProp}(a).  

($ii$): 
In particular,
($i$) implies that $T\!\in\!\mathbb A_\lambda^{(k)}$ for $k\geq h_M(\lambda)$.  
By the definition of $\mathbb A_{\lambda}^{(k)}$,
it thus suffices to show that $\mathbb P_{\lambda^{\to k}} U = 0$ for all $U\neq T$.
This is true by Property~\ref{FilProp}(b). 
\hfill
$\endprf$

As with the definition of the Hall-Littlewood polynomials,
we associate symmetric functions to our super atoms. 
\begin{definition}
With $\digamma$ as in \ref{morphism},
we define the symmetric function atoms by 
\begin{equation}
A_{\lambda}^{(k)}[X;t]=\digamma \left( \mathbb A_{\lambda}^{(k)} \right)\, .
\end{equation}
\label{defatomfunc}
\end{definition}
Properties we have given for the super atoms allow us to deduce
several properties of these functions.  For example,
an immediate consequence of Property~\ref{SAProp}($ii$) is
\begin{property}
When $k$ is large (\,$k \geq h_M(\lambda)$\,), we have
\begin{equation}
A_{\lambda}^{(k)}[X;t] = S_{\lambda}[X] \, .
\end{equation}
\label{corlarge}
\end{property}
\begin{property}
The atoms are linearly independent and have an expansion of the form
\begin{equation}
A_{\lambda}^{(k)}[X;t] = S_{\lambda}[X] + \sum_{\mu > \lambda} 
v_{\mu \lambda}^{(k)}(t) \, S_{\mu}[X] \, ,
\quad\text{where}\;\; v_{\mu \lambda}^{(k)}(t) \in \mathbb N[t]\,.
\label{eqtri}
\end{equation}
\label{proptri}
\end{property}
\noindent {\it Proof.} \quad 
We have shown that $\mathbb A_{\lambda}^{(k)}\subseteq \mathbb H_{\lambda}$.
Thus by Definition~\ref{defatomfunc},
the triangularity of $A_{\lambda}^{(k)}[X;t]$ 
follows from the triangularity of $H_{\lambda}[X;t]$ (see \ref{HLF}).
Further, Property~\ref{SAProp} implies that the
tableau $T$ of shape and evaluation $\lambda$ 
occurs in $\mathbb A_\lambda^{(k)}$
and therefore, $\digamma(T)=S_\lambda[X]$ 
occurs in $A_\lambda^{(k)}$ ($T$ has charge zero). 
\hfill $\endprf$

\section{Main conjecture}

Our work to characterize the atoms was originally 
motivated by the belief that these polynomials play 
an important role in understanding the $q,t$-Kostka coefficients.
More precisely, 
\begin{conjecture}
For any partition $\lambda$ bounded by $k$, 
\begin{equation}
H_{\lambda}[X;q,t]=
\sum_{\mu; \mu_1 \leq k} K_{\mu \lambda}^{(k)}(q,t) \, A_{\mu}^{(k)}[X;t] \, ,
\quad\text{where}\;\;
K_{\mu \lambda}^{(k)}(q,t) \in \mathbb N[q,t] \, .
\label{Mac2atoms}
\end{equation}
\label{conjMac}
\end{conjecture}
\noindent For example, we have
\begin{equation}
\begin{split}
H_{2,1,1}[X;q,t] & = t\, A_{2,2}^{(2)}[X;t]+ (1+q t^2) \, A_{2,1,1}^{(2)}[X;t] + 
q\, A_{1,1,1,1}^{(2)}[X;t] \\
&  = t^2 \, A_{3,1}^{(3)}[X;t]+(t+q t^2) \, A_{2,2}^{(3)}[X;t]+ 
(1+q t^2) A_{2,1,1}^{(3)}[X;t] + q\, A_{1,1,1,1}^{(3)}[X;t] \\
&  = t^3\, A_{4}^{(k\geq 4)}[X;t]+(t+t^2+qt^3) \, A_{3,1}^{(k \geq 4)}[X;t]+(t+q t^2) \, A_{2,2}^{(k \geq 4)}[X;t]\\
& \qquad \qquad+ (1+q t+q t^2) A_{2,1,1}^{(k \geq 4)}[X;t] + 
q\, A_{1,1,1,1}^{(k \geq 4)}[X;t] \, .
\end{split}
\end{equation}
This conjecture implies that the atoms of level $k$ form a basis for $V_k$.  
Further, since the atoms expand positively in terms of Schur functions
\ref{eqtri}, our conjecture also implies Macdonald's positivity conjecture
on the $H_{\lambda}[X;q,t]$ in $V_k$.  
Since Property~\ref{corlarge} gives
\begin{equation}
 K_{\mu \lambda}^{(k)}(q,t)=K_{\mu \lambda}(q,t) \, 
\quad\text{for}\;\; k \geq |\lambda| \, ,
\end{equation}
we see that this conjecture
is a generalization of Macdonald's conjecture. 

In fact, our conjecture refines the original Macdonald conjecture
in the following sense:
substituting \ref{eqtri},
the positive Schur function expansion of atoms,
into \ref{Mac2atoms}, we have
\begin{equation}
\begin{split}
H_{\lambda}[X;q,t] & = \sum_{\mu} K_{\mu \lambda}^{(k)}(q,t) \, 
\sum_{\nu \geq \mu} v_{\nu \mu}^{(k)}(t) \, S_{\nu}[X]\,. 
\end{split}
\end{equation}
On the other hand,  since the $q,t$-Kostka coefficients appear
in the expansion
\begin{equation}
\begin{split}
H_{\lambda}[X;q,t]
& = \sum_{\nu} K_{\nu \lambda}(q,t) \, S_{\nu}[X] \, ,
\end{split}
\end{equation}
we have that
\begin{equation}
\begin{split}
K_{\nu \lambda }(q,t) & = 
\sum_{\mu \leq \nu} K_{\mu \lambda}^{(k)}(q,t) \, v_{\nu \mu}^{(k)}(t) \\
& =  K_{\nu \lambda}^{(k)}(q,t)  +
\sum_{\mu < \nu} K_{\mu \lambda}^{(k)}(q,t) \, v^{(k)}_{\nu \mu}(t) \, .
\end{split} 
\label{pluspet}
\end{equation}
Since $v_{\nu\mu}^{(k)}(t)$ is in $\mathbb N[q,t]$, Conjecture~\ref{conjMac}
implies that
\begin{equation}
 K_{\mu \lambda}^{(k)}(1,1)  \leq K_{\mu \lambda}(1,1) \, , 
\end{equation}
where $K_{\mu \lambda}(1,1)$ is known to be the number of standard tableaux
of shape $\mu$.  
Thus, the problem of finding a combinatorial interpretation for 
the $K_{\mu \lambda}(q,t)$ coefficients (i.e. associating statistics 
to standard tableaux) is reduced to obtaining 
statistics for the fewer $K_{\mu \lambda}^{(k)}(q,t)$.  

Based on our conjecture,
we have the following corollary concerning the
expansion of Hall-Littlewood polynomials in terms of our atoms. 
\begin{corollary}  Assuming Conjecture~\ref{conjMac} holds, we have,
for any partition $\lambda$ bounded by $k$, 
\begin{equation}
H_{\lambda}[X;t]=\sum_{\mu \geq \lambda} K_{\mu \lambda}^{(k)}(t) \, A_{\mu}^{(k)}[X;t] 
\quad
\text{where}\;\;K_{\mu \lambda}^{(k)}(t) \in \mathbb N[t]\,.
\end{equation}
\label{conjHL}
\end{corollary}
\noindent 
If we consider this corollary as the result of 
applying $\digamma$ to an identity on tableaux,
\begin{equation}
\digamma \left(\mathbb H_\lambda\right) 
= \sum_{\mu} K_{\mu\lambda}^{(k)}(t)\,
\digamma \left(\mathbb A_\mu^{(k)}\right) \, ,
\quad \text{where }\; K_{\mu\lambda}(t)\in\mathbb N[t]\, ,
\label{atdecomp}
\end{equation}
then it suggests that the set of 
all tableaux with evaluation 
$\lambda$ can naturally be decomposed 
into subsets that are mapped under $\digamma$
to the atoms $A_\mu^{(k)}[X;t]$.
Here, $K_{\mu \lambda}^{(k)}(1)$ corresponds to the number of 
times such a subset occurs in $\mathbb H_\lambda$ which, by \ref{pluspet},
is such that
\begin{equation}
 K_{\mu \lambda}^{(k)}(1)  \leq K_{\mu \lambda}(1) \,,
\end{equation}
where $K_{\mu \lambda}(1)$ is the number of tableaux
with evaluation $\lambda$ and shape $\mu$.
These subsets  will be called {\it copies} of $\mathbb A_{\mu}^{(k)}$ and
they will provide a natural decomposition for the set of tableaux of a
given evaluation.

\section{Embedded tableaux decomposition}

We expect from \ref{atdecomp} that the set of 
all tableaux with given evaluation can be decomposed 
into subsets associated to our super atoms. 
These subsets will be characterized by 
a cyclage-cocyclage ranked-poset structure \cite{[LS1]}. 

For tableau $T\! = \!xw$ where $x$ is not the smallest letter of $T$,
we define $T'$ to be the unique tableau such that 
$T' \equiv wx$.  The mapping $T \to T'$ is a called a {\it cyclage}
and is such that charge$(T')=$ charge$(T)+1$ if the evaluation of $T$ is
a partition.
For tableau $T\!=\!wx$ where $x$ is not the smallest letter of $T$,
we define $T'$ to be the unique tableau such that $T' \equiv xw$.  
The {\it cocyclage} is the mapping $T \to T'$ 
and is such that charge$(T')=$ charge$(T)-1$ if the evaluation of $T$ is
a partition.

On any collection of tableaux $\mathbb T$ of the same evaluation, we can define
a poset $(\mathbb T, <_{cc})$.  In the case where the evaluation is a partition,
the poset is defined  by linking
any two tableaux $T$ and $T'$ if $T$ is obtained from $T'$,
or vice versa, using either a cyclage or a cocyclage.
In the case where the evaluation is not a partition, we first permute the evaluation
to a partition by using an element $\sigma$ of the symmetric group, and then
 construct the poset  by linking
any two tableaux $T$ and $T'$ if $\sigma T$ is obtained from $\sigma T'$,
or vice versa, using either a cyclage or a cocyclage.
This induces a partial order on $\mathbb T$, 
such that if 
$T\!<_{cc}\!T'$ then charge$(T)\!<$charge$(T')$.
For example, the poset $\left( \mathbb A_{3,2,2,1,1}^{(4)} , <_{cc} \right)$ is
\begin{equation}  
\begin{matrix}
{\rm{charge}} & & & & & & \\ 
3& & & & {\tiny{\tableau*[scY]{ 3 \cr 2 & 2 & 5\cr 1 & 1& 1& 3 & 4}}}
 & & \\
& & & \swarrow  \llap{$ \nearrow $\kern 3pt} & & \nwarrow & \\
2 & & {\tiny{\tableau*[scY]{  4 \cr  3 \cr 2 & 2 & 5\cr 1 & 1& 1 & 3}}} & &        {\tiny{\tableau*[scY]{  3 & 3 \cr 2 & 2 & 5\cr 1 & 1& 1 & 4}}}
 & &
 {\tiny{\tableau*[scY]{ 4 \cr 3 \cr 2 & 2 \cr 1 & 1& 1& 3 & 5}}}\\
 & & & \searrow \llap{$ \nearrow $\kern 3pt}  & \updownarrow & & \updownarrow \\
1 & &   {\tiny{\tableau*[scY]{ 4 \cr 3 & 3 \cr 2 & 2 \cr 1 & 1& 1& 5}}}  & &  {\tiny{\tableau*[scY]{ 4 \cr 3 & 3 \cr 2 & 2 & 5 \cr 1 & 1& 1}}}
 & &
 {\tiny{\tableau*[scY]{ 5 \cr 4 \cr 3 \cr 2 & 2 \cr 1 & 1& 1& 3 }}}\\
& & & \searrow  \llap{$ \nwarrow $\kern 3pt} & & \swarrow & \\
0& & & & {\tiny{\tableau*[scY]{ 5 \cr 4 \cr 3 & 3 \cr 2 & 2 \cr 1 & 1& 1}}} & & 
\end{matrix} \qquad ,
\end{equation}
where the arrows indicate the cyclage and cocyclage relations between tableaux.
\begin{conjecture}
The cyclage and cocyclage induce a \underline{connected}
ranked-poset structure on the set of tableaux contained
in a given super atom $\mathbb A_{\lambda}^{(k)}$.
\label{conjposet}
\end{conjecture}

Given a collection of tableaux $\mathbb T$, the Hasse diagram of the poset 
$(\mathbb T,<_{cc})$ with vertices labeled by shapes of the corresponding 
tableaux will be denoted $\Gamma_{\mathbb T}$.  We use the symbol
$\Gamma_\lambda^{(k)}$ when $\mathbb T = \mathbb A_\lambda^{(k)}$.
For example, the Hasse diagram $\Gamma_{3,2,2,1,1}^{(4)}$, 
associated to $\mathbb A_{3,2,2,1,1}^{(4)}$, is
\begin{equation}  
\begin{matrix}
{\rm{charge}} & & & & & & \\ 
c+3& & & & {\tiny{\tableau*[scY]{  \cr  &  & \cr  & & &  & }}}
 & & \\
& & & \diagup   & & \diagdown & \\
c+2 & & {\tiny{\tableau*[scY]{   \cr   \cr  &  & \cr  & &  & }}} & &        {\tiny{\tableau*[scY]{   &  \cr  &  & \cr  & &  & }}}
 & &
 {\tiny{\tableau*[scY]{  \cr  \cr  &  \cr  & & &  & }}}\\
 & & & \diagdown \llap{$ \diagup $\kern 0pt} & | & & | \\ 
c+1 & &   {\tiny{\tableau*[scY]{  \cr  &  \cr  &  \cr  & & & }}}  & &  
{\tiny{\tableau*[scY]{  \cr  &  \cr  &  &  \cr  & & }}}
 & &
 {\tiny{\tableau*[scY]{  \cr  \cr  \cr  &  \cr  & & &  }}}\\
& & & \diagdown & & \diagup & \\
c& & & & {\tiny{\tableau*[scY]{  \cr  \cr  &  \cr  &  \cr  & & }}} & & 
\end{matrix} \qquad 
\end{equation}
We can now define the subsets associated to our atoms. 
\begin{definition}
If a set of tableaux $\mathbb T$  
has the properties:
\begin{equation}
\begin{split}
\quad &1.\;T \text{ is the tableau of minimal charge in }
\mathbb T 
\qquad \qquad \qquad \qquad
\qquad \qquad
\qquad 
\\
&
2.\; \Gamma_{\mathbb T} = \Gamma_{\text{shape}(T)}^{(k)}\,,
\qquad \qquad \qquad \qquad
\qquad \qquad 
\qquad 
\end{split}
\label{copy1}
\end{equation}
{\vskip-0.2in
\noindent then this set is called a copy of the atom 
$\mathbb A_{\text{shape}(T)}^{(k)}$
and is denoted $\mathbb A_T^{(k)}$.
}
\end{definition}
\noindent
This given, if the posets are connected 
(Conjecture~\ref{conjposet}) then the charges 
associated to the elements of a super atom 
$\mathbb A_{\lambda}^{(k)}$ differ from those 
of a copy atom $\mathbb A_{T}^{(k)}$
by a common factor.  
Furthermore, since there is a unique element of
zero charge in $\mathbb A_{\lambda}^{(k)}$
(the tableau with shape and evaluation $\lambda$),
then it is the minimal element in $\mathbb A_{T}^{(k)}$
and we have
\begin{equation}
\digamma\left(\mathbb A_T^{(k)}\right) = 
t^{\text{charge}(T)}
A_\lambda^{(k)}[X;t] \, , \quad
\text{where}\;\; \text{shape}(T)=\lambda\, .
\label{digamat}
\end{equation}
Note also that the tableaux in 
$\mathbb A_{\lambda}^{(k)}$ have evaluation 
$\lambda$ while those in 
$\mathbb A^{(k)}_T$ have evaluation given by the
evaluation of $T$. That is,
$\mathbb A_{\lambda}^{(k)}=\mathbb A_{T}^{(k)}$
only if $T$ is of shape and evaluation $\lambda$.
The copies of $\mathbb A_{3,2,2,1,1}^{(4)}$ 
include $\mathbb A_{863925147}^{(4)}$, given by
\begin{equation}
\begin{matrix}
{\rm{charge}} & & & & & & \\ 
13& & & & {\tiny{\tableau*[scY]{  3 \cr 2 & 5 & 7 \cr 1 & 4& 6 & 8 & 9}}} & & \\
& & & \swarrow \llap{$ \nearrow $\kern 3pt} & & \nwarrow & \\
12 & & {\tiny{\tableau*[scY]{ 9 \cr 3  \cr 2 & 5 & 7 \cr 1 & 4& 6& 8}}} & & {\tiny{\tableau*[scY]{  3 & 9 \cr 2 & 5 & 7 \cr 1 & 4& 6 & 8}}} & &
{\tiny{\tableau*[scY]{  6 \cr 3 \cr 2 & 5 \cr 1 & 4& 7& 8 & 9}}} \\
 & & & \searrow  \llap{$ \nearrow $\kern 3pt} & \updownarrow & & \updownarrow \\
11 & & {\tiny{\tableau*[scY]{  6 \cr  3 & 9 \cr 2 & 5 \cr 1 & 4& 7 & 8}}} & & 
{\tiny{\tableau*[scY]{  8 \cr 3 & 9 \cr 2 & 5 & 7\cr 1 & 4& 6 }}} & &
{\tiny{\tableau*[scY]{ 9 \cr 6 \cr 3 \cr 2 & 5 \cr 1 & 4& 7& 8 }}} \\
& & & \searrow  \llap{$ \nwarrow $\kern 3pt} & & \swarrow & \\
10& & & & {\tiny{\tableau*[scY]{ 8 \cr 6 \cr 3 & 9 \cr 2 & 5 \cr 1 & 4& 7}}} & &
\end{matrix} 
\end{equation} 

It appears that there is a unique way to decompose 
the set of all tableaux $\mathbb H_{\mu}$ 
into atoms of level $k \geq \mu_1$. 
More precisely, we let ${\cal C}_\mu^{(k)}$
denote the collection of all tableaux $T$ 
with evaluation $\mu$ where $\mathbb A_T^{(k)}$ is a copy of a super atom of level $k$.  
Then

\begin{conjecture}
For any partition $\mu$ bounded by $k$, we have
\begin{equation}
\mathbb H_\mu
=\sum_{T\in{\cal C}_\mu^{(k)}}
\mathbb A_T^{(k)}
\,.
\label{HinAt}
\end{equation}
\end{conjecture}
\noindent From Corollary~\ref{conjHL} and \ref{digamat}, an implication of this identity
under the mapping $\digamma$ is:
\begin{corollary}
The $k$-Kostka-Foulkes polynomials are
simply
\begin{equation}
K_{\lambda\mu}^{(k)}(t) 
=\sum_{
\begin{subarray}{c}
T \in {\cal C}_\mu^{(k)}\\
\text{shape}(T)=\lambda
\end{subarray}
}
t^{\text{charge}(T)}
\,.
\label{kostkafoulkes}
\end{equation}
\end{corollary}

One method to obtain the set
${\cal C}_\mu^{(k)}$ is as follows:
The element of minimal charge in $\mathbb H_\mu$ 
has shape $\mu$ and is thus also the minimal element
of $\mathbb A_{\mu}^{(k)}$.  
Remove from $\mathbb H_\mu$, all tableaux in $\mathbb A_{\mu}^{(k)}$.
Choose a tableau $T$ with minimal charge from those that remain.
$T$ must index a copy of the atom $\mathbb A_{\lambda}^{(k)}$
where $\lambda=\text{shape}(T)$.
From the Hasse diagram $\Gamma_{\lambda}^{(k)}$, 
it is possible to find and remove all tableaux in 
the atom $\mathbb A_{T}^{(k)}$. 
Repeat this procedure always on an element of minimal charge 
in the resulting sets.  The collection of these minimal elements
is ${\cal C}_\mu^{(k)}$.

Evidence suggests that this method also
provides a direct decomposition of any copy atom of a level $k$ 
into copy atoms of level $k'>k$. 
\begin{conjecture}  For any  atom $\mathbb A_T^{(k)}$ 
such that $\shape(T)$ is bounded by $k$, and any $k' > k$,
\begin{equation}
\mathbb A_T^{(k)} = \sum_{T' \in {\cal D}^{(k \to k')}_{T}} \mathbb A_{T'}^{(k')} \, ,
\end{equation}
for some collection of tableaux ${\cal D}^{(k \to k')}_{T}$.
\label{conjdecomp}
\end{conjecture}
On the level of functions, this translates into
\begin{corollary} If $\lambda$ is a partition bounded by $k$, and $k' > k $,  then
\begin{equation}
A_{\lambda}^{(k)}[X;t] = A_{\lambda}^{(k')}[X;t]+
\sum_{\mu > \lambda} v_{\mu \lambda}^{(k \to k')}(t) \, A_{\mu}^{(k')}[X;t] \, ,
\quad
\text{where}\;\; v_{\mu \lambda}^{(k \to k')}(t) \in \mathbb N[t]\,.
\end{equation}
\end{corollary}
\noindent This conjecture is a generalization of the result
presented in Property~\ref{proptri} since
we recover $v_{\mu \lambda}^{(k \to k')}(t)=v_{\mu \lambda}(t)$
when $k' \geq |\lambda|$.

\begin{figure}[htb]
\begin{center}
\epsfig{file=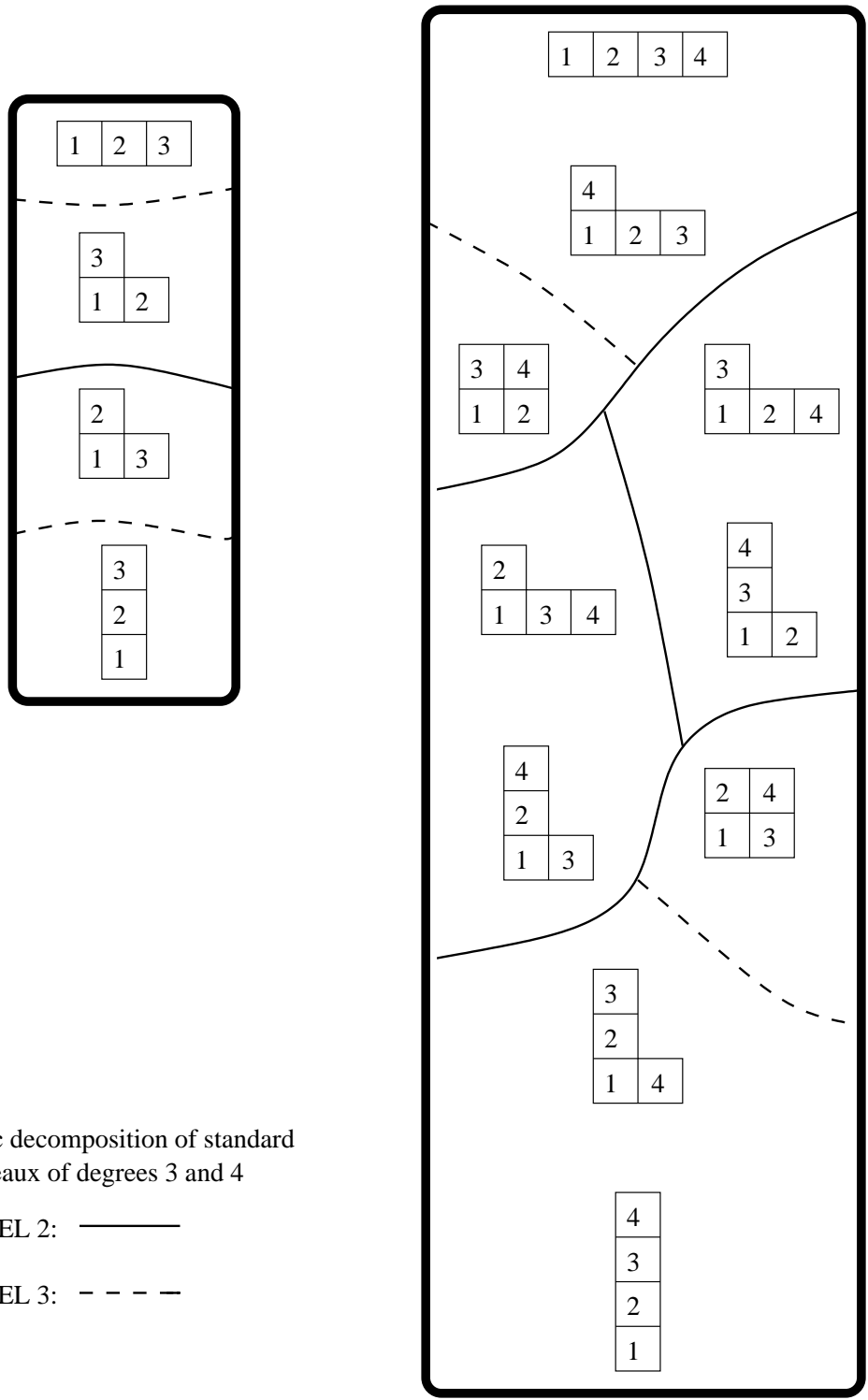}
\end{center}
\caption{In order to read the decomposition of a given level $k$, one must consider the lines associated to all the levels which are not bigger than $k$.  That is, when doing the decomposition from one level to the other, lines
are added without ever being removed.  Thus for instance the tableau $2413$ and $3214$ are in the same atom up to level 2, and in different
atoms for any higher levels. 
}
\label{case34}
\end{figure}

\begin{figure}[htb] 
\begin{center}
\epsfig{file=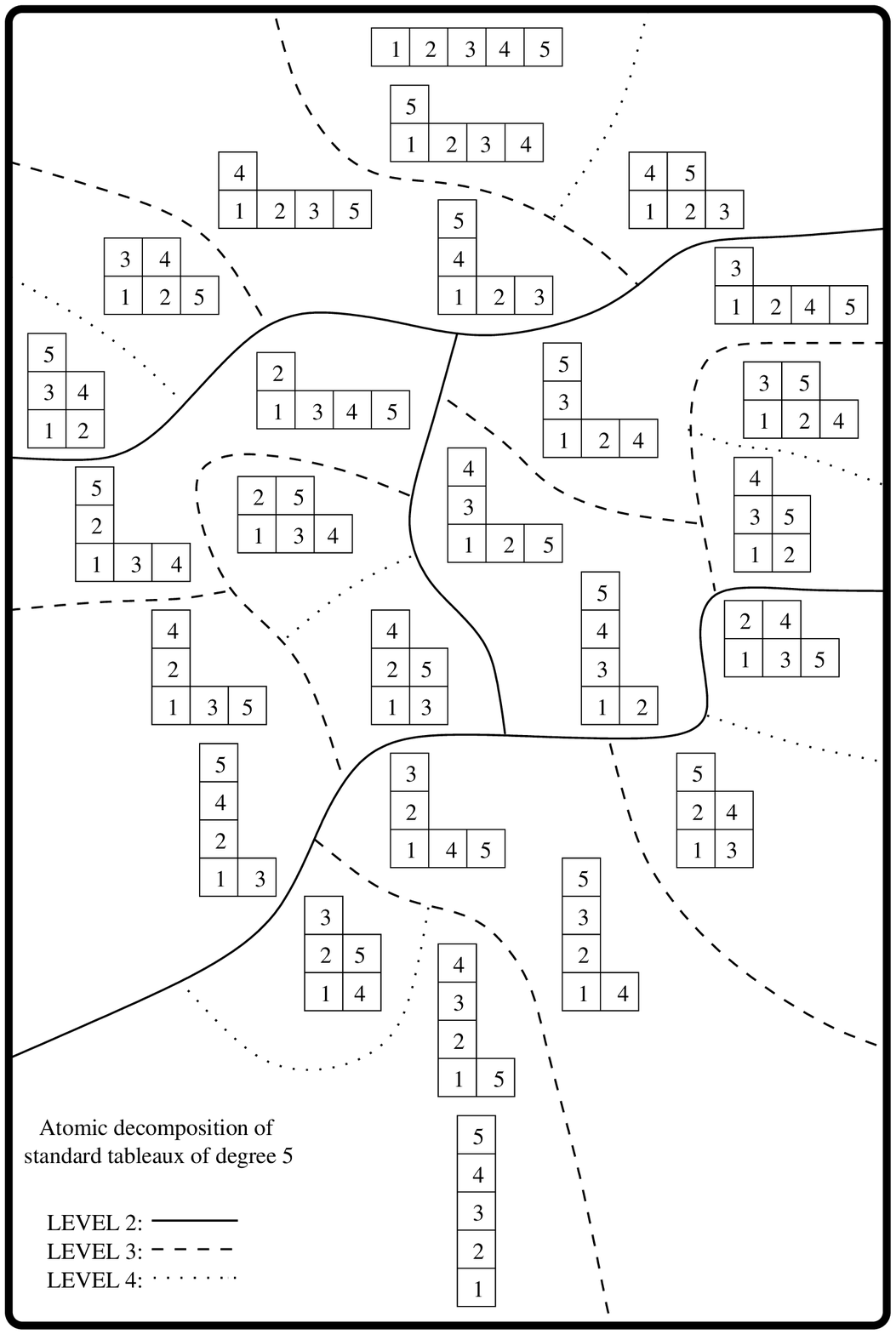}
\end{center}
\caption{See Figure~\ref{case34} for details on how to read the figure.}
\label{case5}
\end{figure}

For examples that support the preceding conjectures, 
refer to Figures~\ref{case34} and \ref{case5}.
These figures also suggest that the number of elements in an atom, 
at increasing charges, forms a unimodal sequence.  Since an atom has a unique minimal element, these sequences always start with 1.
\begin{conjecture}  Given any atom
\begin{equation}
A_{\lambda}^{(k)}[X;t] = \sum_{\mu \geq \lambda} v_{\mu \lambda}^{(k)}(t) 
\, S_{\mu}[X] \, ,
\end{equation}
the numbers
\begin{equation}
\#_i = \sum_{\mu \geq \lambda} v_{\mu \lambda}^{(k)}(t) \Big |_{t^i} \, ,
\end{equation}
are such that $[\#_0,\#_1,\dots]$ is a unimodal sequence.
\end{conjecture}
For example, the unimodal sequence associated to
$A_{3,2,2,1,1,1}^{(4)}[X;t]$
is $[1,3,5,5,3,1]$:
\begin{equation}
\begin{split}
A_{3,2,2,1,1,1}^{(4)}[X;t] \, & = \,  
S_{{3,2,2,1,1,1}} +
t \, S_{{4,2,1,1,1,1}}
+ t \, S_{{3,3,2,1,1}}
+ \left (t+{t}^{2}\right )S_{{4,2,2,1,1}}
+t^2\,  S_{{3,3,3,1}} \\
& +t^2 \, S_{{4,3,1,1,1}}   
+ \left ({t}^{2}+{t}^{3}\right )S_{{5,2,1,1,1}}
+\left ({t}^{2}+{t}^{3}\right )S_{{4,3,2,1}}
+ t^3 \, S_{{5,2,2,1}} \\
&+t^3 \, S_{{4,3,3}}
+\left ({t}^{3}+{t}^{4} \right )S_{{5,3,1,1}} 
+t^4 \, S_{{6,2,1,1}} 
+t^4 \, S_{{5,3,2}}  
+t^5 \, S_{{6,3,1}}
\, .
\end{split}
\end{equation}
We will see later (Corollary~\ref{cormin}) that these sequences
also end with a one,
that is, an atom has a unique element of maximal charge.
We will also provide a way to obtain the shape of this maximal element.
Note that the sequences are not necessarily symmetric.  For instance,
from Figure~\ref{case5} we see that
the sequence associated to $A_{1,1,1,1,1}^{(2)}[X;t]$ is [1,1,2,2,1].

We finish this section by stating a conjecture 
that reiterates the importance of the atoms as
a natural basis for $V_k$.
\begin{conjecture}
For any two alphabets $X$ and $Y$, 
\begin{equation}
A_\lambda^{(k)}[X+Y;t]=\sum_{|\mu|+|\rho|=|\lambda|}
g_{\mu \rho}^\lambda(t)\,A_{\mu}^{(k)}[X;t]\,A_{\rho}^{(k)}[Y;t] \, ,
\end{equation}
with $g_{\mu\rho}^\lambda(t)\in \mathbb N[t]$. 
\end{conjecture}
It is important to note that the positivity of the coefficients $g_{\mu\rho}^{\lambda}(t)$
appearing here is a natural property of Schur functions that is not shared by the 
Hall-Littlewood or Macdonald functions. 

\section{Irreducible atoms}

We have now seen that the super atoms can be constructed by
generating sets of tableaux with promotion operators $\mathbb B_{r}$
and then eliminating undesirable elements using
the projection operators $\mathbb P_{\lambda^{\to k}}$. 
Further, we have given a method to obtain copies of the super atoms 
allowing us to decompose the set of all tableaux with a given 
evaluation and to provide 
natural properties on the functions $A_\lambda^{(k)}[X;t]$.

Remarkably, it appears that there is a method to construct many of 
the atoms without generating any undesirable elements.
In fact, what could be seen as the `DNA' of our atoms 
is a subset of irreducible atoms for each $V_k$, from 
which all successive atoms of $V_k$ may be obtained by 
simply applying a generalized version 
of the promotion operators.

To be more precise, let a rectangular partition of the form
$(\ell^{k+1-\ell})$ be referred to as a $k$-rectangle
and a partition with no more than $i$ parts equal to $k-i$  
be called $k$-irreducible.

\begin{definition}
The collection of $k$-irreducible atoms
is composed of atoms indexed by $k$-irreducible partitions.
If an atom is not irreducible, then it is said to be reducible.
\label{defindecom}
\end{definition}
\begin{property}  
There are $k!$ distinct $k$-irreducible partitions.
\label{lemmakfac}
\end{property}

\noindent {\it Proof.} \quad 
A partition $\lambda$ is $k$-irreducible if and only if 
$\lambda$ has no more than $i$ parts equal to $k-i$.
There are obviously $k!$ such partitions.  \hfill $\endprf$

\noindent The irreducible atoms of level 1,2 and 3 are
\begin{equation}
\begin{split}
k=1 \, &:  \qquad A_0^{(1)} \, ,\\
k=2 \, &: \qquad A_0^{(2)} \, , \quad A_1^{(2)} \, ,\\
k=3 \, &: \qquad A_0^{(3)} \, , \quad A_1^{(3)}\, ,
\quad A_{2}^{(3)}\, , \quad A_{1,1}^{(3)} \, , \quad A_{2,1}^{(3)} \, , \quad A_{2,1,1}^{(3)}\, .
\end{split}
\end{equation}

Any $k$-bounded partition $\mu$ is of the form
$\mu = \lambda \cup R_1 \cup \cdots \cup R_n$, where $\lambda$ is a
$k$-irreducible  partition and
$R_1,\dots, R_n$ is a sequence of $k$-rectangles.
In fact, any atom of level $k$ can be obtained from a
$k$-irreducible atom by the application 
of certain generalized promotion operators
that are indexed by $k$-rectangles.

Before we can introduce these promotion operators, 
we need to define an operation which generalizes $\sigma_i$.
We define $\sigma_i^{(h)}$ to send a word $w$ of evaluation 
$(\rho_1,\rho_2,\dots)$ to a word $w'$ of evaluation
$(\rho_1,\dots,\rho_{i-1},\rho_{i+1},\dots,\rho_{i+h},\rho_i,\rho_{i+h+1},\dots)$.
The operation $\sigma_i^{(h)}$ only acts on the subword $w_{\{i,\dots,i+h\}}$, and can
thus be defined in generality from the special case $i=1$. 
Let $w$ be a word in $1,\dots,h+1$ and
let $\bigl(P(w),Q(w)\bigr)$
denote the pair of tableaux in RS-correspondence with $w$ (see \ref{RS}).
If $w''$ is the word obtained from $w$ by first erasing all occurrences of 
the letter 1 and then decreasing the remaining letters by 1,
then the shape of $P(w'')$ differs from
that of $P(w)$ by a horizontal strip. 
Let $T'$ be the tableau obtained from
$P(w'')$ by filling the horizontal strip with $(h+1)$'s, and let 
$w'$ be the word which is in RS-correspondence with the pairs 
of tableaux $\bigl(P(w'),Q(w')\bigr)=\bigl(T',Q(w)\bigr)$.  
We now define $\sigma_1^{(h)}$ by
setting
\begin{equation}
\sigma_1^{(h)} ( w ) = w'\, .
\end{equation}
It can  be shown that $\sigma_i^{(1)}=\sigma_i$,
and thus  $\sigma_i^{(h)}$ generalizes $\sigma_i$.
Note that
$\sigma_i^{(h)}$ happens to be a special case of an operation
defined in \cite{[S2]}.

The rectangular promotion operators are defined in a manner similar to 
the promotion operators.
That is, on a tableau $T$ of evaluation $(\lambda_1,\ldots,\lambda_m)$,
\begin{equation}
\mathbb B_{(\ell^h)} ( T )= 
\sigma_1^{(h)}\cdots\sigma_m^{(h)} \, \mathbb R_{(\ell^h)}
\, T\, 
\label{RecOp}
\end{equation}
generates a sum of tableaux with evaluation
$(\ell^h,\lambda_1,\ldots,\lambda_m)$
by applying a rectangular analogue of $\mathbb R_r$. 
This operator, $\mathbb R_{(\ell^h)}$, acts by adding to $T$, 
a horizontal $\ell$-strip of the letter
$m+1$, a horizontal $\ell$-strip of $m+2$,\ldots,
and a horizontal $\ell$-strip of $m+h$ in 
all possible ways such that the tableaux are 
Yamanouchi in the added letters
\footnote{
This is the multiplication involved in
computing the Littlewood-Richardson coefficients in
the product of a Schur function of the shape of $T$ by a Schur function
indexed by a rectangular partition.}.
Since $\sigma_i^{(1)}=\sigma_i$ for $h=1$,
we recover the previously defined
promotion operator $\mathbb B_{\ell}$.
\begin{equation}
\begin{split}
\mathbb B_{(2^3)} \, \, {\tiny{\tableau*[scY]{ 2 \cr 1 & 2}}}
 & \, = \,     \sigma_{1}^{(3)} \, \sigma_2^{(3)} \,
\mathbb R_{(2^3)} \, \, {\tiny{\tableau*[scY]{ 2 \cr 1 & 2}}}
\,  = \,   \sigma_{1}^{(3)} \sigma_2^{(3)} \,
\left(\,  {\tiny{\tableau*[scY]{ 5 & 5 \cr 2 & 4 & 4 \cr 1 & 2 & 3 & 3 }}}
+{\tiny{\tableau*[scY]{ 5 \cr 4 & 5 \cr 2 & 4 \cr 1& 2 & 3 & 3}}}
+{\tiny{\tableau*[scY]{ 5 & 5 \cr 4 & 4 \cr 2 & 3 \cr 1 &2 & 3}}}
+{\tiny{\tableau*[scY]{ 5 \cr 4 & 5 \cr 2 & 3 & 4\cr 1 &2 & 3}}}
+{\tiny{\tableau*[scY]{ 5 \cr 4 \cr 3 & 5 \cr 2 & 4 \cr 1 &2 & 3}}}
+{\tiny{\tableau*[scY]{ 5 \cr 4 & 5 \cr 3 & 4 \cr 2 & 3 \cr 1 & 2}}} \, \right)\\
& = \quad {\tiny{\tableau*[scY]{ 3& 3 \cr 2 & 2 & 5 \cr 1 & 1 & 4 & 5 }}}+
{\tiny{\tableau*[scY]{ 4 \cr 3 & 3 \cr 2 & 2 \cr 1& 1 & 5 & 5}}}+
{\tiny{\tableau*[scY]{ 4 & 5 \cr 3 & 3 \cr 2 & 2 \cr 1 &1 & 5}}}+
{\tiny{\tableau*[scY]{ 5 \cr 3 & 3 \cr 2 & 2 & 5\cr 1 &1 & 4}}}+
{\tiny{\tableau*[scY]{ 5 \cr 4 \cr 3 & 3 \cr 2 & 2 \cr 1 &1 & 5}}}+
{\tiny{\tableau*[scY]{ 5 \cr 4 & 5 \cr 3 & 3 \cr 2 & 2 \cr 1 & 1}}} \, .
\end{split}
\label{BonA2}
\end{equation} 
In fact, it seems that the inverse of 
the $\mathbb B_{(\ell^h)}$
action on an atom of level $k$ is 
simply rectangular-katabolism,  
\begin{conjecture}  
If $\tau$ is a translation of the letters in $\mathbb A_T^{(k)}$,
we have
\begin{equation}
\mathbb K_{(\ell^{k-\ell+1})}\,\mathbb B_{(\ell^{k-\ell+1})} \, 
\mathbb A_T^{(k)}
\, = \,
\tau \mathbb A_T^{(k)}
\, .
\end{equation}
\label{conjkata}
\end{conjecture}
\noindent The tableaux in \ref{BonA2} are sent, under $\mathbb K_{(2^3)}$,
to $\,{\tiny{\tableau*[scY]{  5 \cr 4& 5}}}\,$
\big(a translation of the atom
$\mathbb A_{212}^{(3)}\,$ \big).  Note that in general 
$\mathbb K_{(\ell^h)} \, \mathbb B_{(\ell^h)} \, T$ is not necessarily equal
to a translation of the letters in $T$.

This conjecture supports the very important 
idea that any atom can be obtained from an 
irreducible atom simply by applying a sequence of
rectangular promotion operators.   That is,
\begin{conjecture} 
The operator $\mathbb B_{(\ell^{k-\ell+1})}$ acts 
on any copy $\mathbb A_T^{(k)}$ of $\mathbb A_{\lambda}^{(k)}$
by
\begin{equation}
\mathbb B_{(\ell^{k-\ell+1})} \, \mathbb A_T^{(k)} = \mathbb A_{T'}^{(k)} \, , 
\end{equation}
for a tableau $T'$ of shape $\lambda \cup (\ell^{k-\ell+1})$. 
\label{conjrec}
\end{conjecture}
\noindent For instance, by applying $\mathbb B_{(3)}$ 
to $\mathbb A_{213}^{(3)}= {\tiny{\tableau*[scY]{  2 \cr 1 & 3}}} \, $,
we obtain a copy of $\mathbb A_{3,2,1}^{(3)}$:
\begin{equation}
\mathbb B_{(3)} \, {\tiny{\tableau*[scY]{  2 \cr 1 & 3}}} 
=
{\tiny{\tableau*[scY]{  3 \cr 2 & 4 \cr 1 & 1 & 1}}}+
{\tiny{\tableau*[scY]{  2 \cr 1 & 1 & 1 & 3 & 4}}}+
{\tiny{\tableau*[scY]{  4 \cr 2 \cr 1 & 1 & 1 & 3}}}+
{\tiny{\tableau*[scY]{  2 & 4\cr 1 & 1 & 1 &3}}} 
=\mathbb A_{324111}^{(3)} \,.
\label{exBonA}
\end{equation}

Conjecture~\ref{conjrec} not only reveals the importance of the set
of irreducibles, but also provides a convenient way to obtain   
copies using a simple transformation on tableau.
Given a tableau, the transformation $\mathbb L_{T}$ is defined 
by replacing the $\text{shape}(T)$-subtableau with $T$ and 
then adjusting the remaining entries to start with $s+1$, 
where $s$ is the largest letter of $T$.
$\mathbb L_T$ satisfies several properties on the set
of tableaux, denoted $\mathbb H_{\mu|_{(\ell^h)}}$, with evaluation $\mu$
(not necessarily a partition) whose restriction to the $h$ smallest letters
gives exactly the subtableau of shape and evaluation
$(\ell^h)$.

\begin{property}
Let $\mathbb T\subseteq\mathbb H_{\mu|_{(\ell^h)}}$
be a set containing a unique element of 
minimal charge and whose poset $(\mathbb T,<_{cc})$ 
is connected.  For any tableau $T$ of shape $(\ell^h)$,  
$\overline{\mathbb T} = \mathbb L_{T} \mathbb T$
satisfies

\noindent\;\; 1. $\Gamma_{\mathbb T} = \Gamma_{\overline{\mathbb T}}$. 

\noindent\;\; 2. If $U$ is the element of minimal charge in $\mathbb T$,
then $\mathbb L_T U$ is the minimal element in $\overline{\mathbb T}$.
\label{theorLop}
\end{property}
\noindent In particular, if we assume  that 
$\mathbb A_{T'}^{(k)}=\mathbb B_{(\ell^{k-\ell+1})}\mathbb A_T^{(k)}$,
then $\mathbb A_{T'}^{(k)}\subseteq \mathbb H_{\mu|_{(\ell^h)}}$
and by Conjecture~\ref{conjposet}, 
the poset $(\mathbb A_{T'}^{(k)},<_{cc})$ is connected.  
Therefore, 
for any tableau $U$ of shape $(\ell^{k-\ell+1})$, 
$\mathbb L_{U} \mathbb A_{T'}^{(k)}$ 
satisfies the conditions above.
However, these are exactly the conditions 
on a copy of an atom (see \ref{copy1})
and thus 
$\mathbb L_{U} \mathbb A_{T'}^{(k)}$ 
is a copy.
\begin{corollary}
If $\mathbb A_{T_1}^{(k)}=\mathbb B_{(\ell^h)}\mathbb A_{T_2}^{(k)}$ 
is a copy of $\mathbb A_{\lambda}^{(k)}$,
then for each tableau $T$ of shape $(\ell^h)$,
\begin{equation}
\mathbb A_{\mathbb L_{T} T_1}^{(k)}
\,=\,
\mathbb L_{T}\, \mathbb A_{T_1}^{(k)} 
\end{equation}
is another copy of $\mathbb A_{\lambda}^{(k)}$.
\end{corollary}

In example \ref{exBonA}, we let $111\to 123$ to obtain 
another copy of $\mathbb A_{3,2,1}^{(3)}$: 
\begin{equation}
\mathbb A_{546123}^{(3)}
\, =\,
{\tiny{\tableau*[scY]{  5 \cr 4 & 6 \cr 1 & 2 & 3}}}+
{\tiny{\tableau*[scY]{  4 \cr 1 & 2 & 3 & 5 & 6}}}+
{\tiny{\tableau*[scY]{  6 \cr 4 \cr 1 & 2 & 3 & 5}}}+
{\tiny{\tableau*[scY]{  4 & 6\cr 1 & 2 & 3 &5}}} \, .
\label{brecex}
\end{equation}

\noindent {\it Proof of Property~\ref{theorLop}.} \quad 
Let $T_{(\ell^h)}$ be the tableau of shape and evaluation $(\ell^h)$.
Every element $U\in\mathbb T$ contains $T_{(\ell^h)}$, 
and any letter in $U/T_{(\ell^h)}$ is larger than
those in $T_{(\ell^h)}$. 
Therefore the cyclage or cocyclage that links two elements 
$U$ and $U'$ of $\mathbb T$ does not involve the letters 
in $T_{(\ell^h)}$ and we can thus change the content of this 
subtableau without affecting the cyclage-cocyclage relations, 
as long as the new subtableau also contains the smallest letters.  
Hence, the Hasse diagrams of the posets 
$(\mathbb T,<_{cc})$ and $(\overline {\mathbb T}, <_{cc})$
are identical.  
The second condition follows from the 
connectedness of the poset $(\mathbb T,<_{cc})$
which implies that cyclage-cocyclage relations 
could not be preserved if the element of minimal charge in $\overline{\mathbb T}$
was not $\mathbb L_T$ applied on the element of minimal charge of $\mathbb T$.
 \hfill $\endprf$

Another consequence of Conjecture~\ref{conjrec}
arises from the case $t=1$.
Here, the action of $\mathbb B_{(\ell^h)}$ 
on a tableau $T$ is associated to the multiplication of
the Schur functions $S_{(\ell^h)}S_{\text{shape}(T)}$.
\begin{corollary} 
Assuming Conjecture~\ref{conjrec} holds, if we 
let $A_\lambda^{(k)}=A_\lambda^{(k)}[X;1]$,
then 
\begin{equation}
S_{(\ell^{k-\ell+1})} \, A_{\lambda}^{(k)} = 
A_{\lambda \cup (\ell^{k-\ell+1})}^{(k)} \, .
\end{equation}
\label{conjrecschur}
\end{corollary}

We have now seen that any atom can be 
understood as the application of
rectangular promotion operators
to an irreducible component.  
Our study is thus reduced to 
examining the irreducibles
(atoms of level $k$ that cannot be obtained 
by applying $k$-rectangular operators
to a smaller atom).  
Interestingly, we can obtain the
level $k$ atom indexed by the irreducible partition 
of maximal degree,
\begin{equation}
\lambda_M = \bigl( (k-1)^{1},(k-2)^{2},\cdots,1^{{k-1}} \bigr)
\, ,
\end{equation}
by a recursive application of 
$(k-1)$-rectangular promotion operators 
on the empty tableau.
\begin{conjecture}
The maximal irreducible atom of level $k$
is an atom of level $k-1$;
\begin{equation}
A_{\lambda_M}^{(k)}[X;t] = A_{\lambda_M}^{(k-1)}[X;t] \, .
\end{equation}
Furthermore, from Conjecture~\ref{conjrec}, 
this atom is simply
\begin{equation}
A_{\lambda_M}^{(k)}[X;t] = A_{\lambda_M}^{(k-1)}[X;t] = 
\digamma \left( \mathbb B_{(k-1)} \mathbb B_{((k-2)^2)} \cdots 
\mathbb B_{(1^{k-1})} \, \mathbb H_0 \right) \, ,
\end{equation}
\end{conjecture}
\noindent
For example, the atom $A_{2,1,1}^{(3)}[X;t]$ is given by 
\begin{equation}
A_{2,1,1}^{(3)}[X;t]  = \digamma \left( \mathbb B_{(2)} \,  
\mathbb B_{(1^2)} \, \mathbb H_0\right) \,  
  = \digamma \left(\, {\tiny{\tableau*[scY]{  3 \cr 2 \cr 1 & 1}}} +
{\tiny{\tableau*[scY]{  2 \cr 1 & 1 & 3}}} \,   \right) \,  = 
S_{2,1,1}[X;t]+t \, S_{3,1}[X;t] \, .
\end{equation}

When $t=1$, $V_k=\{H_\lambda[X;t]\}_{\lambda_1\leq k}$ reduces to 
the polynomial ring $\mathbb Q[h_1,\dots,h_k]=V_k(1)$. 
If $\mathcal I_k$ denotes the ideal generated by the 
$k$-rectangular Schur functions $S_{(\ell^{k+1-\ell})}$, we have the following 
proposition.
\begin{proposition}  The homogeneous functions indexed by $k$-irreducible partitions
form a basis of the quotient ring $V_k(1)/\mathcal I_k$. 
\end{proposition}
\noindent {\it Proof.} \quad For a partition $\lambda$ bounded by $k$, 
we set
\begin{equation}
\tilde h_{\lambda} = 
\begin{cases}
h_{\lambda} & \text{if $\lambda$ is $k$-irreducible} \\
\tilde h_{\mu+(\ell^{k+1-\ell})}=S_{(\ell^{k+1-\ell})} \, \tilde h_{\mu} &  
\text{if $\lambda=\mu \cup (\ell^{k+1-\ell})$} 
\end{cases} \,. 
\end{equation}
These elements are indexed by $k$-bounded partitions
and thus, if independent, span a space with the same dimension as $V_k(1)$.  
In fact, the $\tilde h_\lambda$ form a basis for $V_k(1)$ since 
$S_{(\ell^{k+1-\ell})}={\rm {det}}\left( h_{\ell-i+j} 
\right)_{1 \leq i,j \leq k+1-\ell}$ implies that
$\tilde h_{\lambda} \in V_k(1)$; and 
$S_{\lambda}=h_{\lambda}+\sum_{\mu> \lambda} c_{\mu \lambda} h_{\mu}$
gives $\tilde h_{\lambda}=h_{\lambda}+\sum_{\mu > \lambda} d_{\mu \lambda} h_{\mu}$, which implies that they are independent.  

First note that the $\tilde h_{\lambda}$ span the quotient ring $V_k(1)/\mathcal I_k$ because they span $V_k(1)$. Since by definition $\tilde h_{\mu} \equiv 0$ in the quotient ring
when $\mu$ is not $k$-irreducible,   
the $\tilde h_\lambda$ indexed by $k$-irreducible partitions will
form a basis for the quotient ring $V_k(1)/\mathcal I_k$ if they are 
independent in $V_k(1)/\mathcal I_k$.  
Let $\mathcal S$ be the set of all $k$-irreducible partitions.  
If, in $V_k(1)/\mathcal I_k$,  
we have
\begin{equation}
\sum_{\lambda \in \mathcal S} d_{\lambda} \, \tilde h_{\lambda}= 0
\, ,
\label{zero}
\end{equation} 
then, in $V_k(1)$, we must have
\begin{equation}
\sum_{\lambda \in \mathcal S} d_{\lambda} \, \tilde h_{\lambda}= \sum_i {C_i} 
\, S_{(i^{k+1-i})}
\, ,
\end{equation} 
for some $C_i\in V_k(1)$. 
Further, since $C_i=\sum_{\mu} c_{i,\mu} \tilde h_{\mu}$ for some $c_{i,\mu}$,
we have 
\begin{equation}
\sum_{\lambda \in \mathcal S} d_{\lambda} \, \tilde h_{\lambda} = 
\sum_{i,\mu} c_{i,\mu} \tilde h_{\mu} S_{(i^{k+1-i})} 
                   = \sum_{i,\mu} c_{i,\mu} \tilde h_{\mu+(i^{k+1-i})} 
\, .
\label{quotient}
\end{equation}
The basis elements appearing in the l.h.s of \ref{quotient} are each
indexed by $k$-irreducible partitions whereas 
those appearing in the r.h.s are indexed by non-$k$-irreducible partitions.  
Therefore, $d_{\lambda}=0$ for all $\lambda$ 
and by \ref{zero},
this proves that the $\tilde h_{\lambda}$ indexed by $k$-irreducible partitions
are independent in $V_k(1)/\mathcal I_k$.
\hfill $\endprf$

\begin{figure}[htb]
\begin{center}
\input{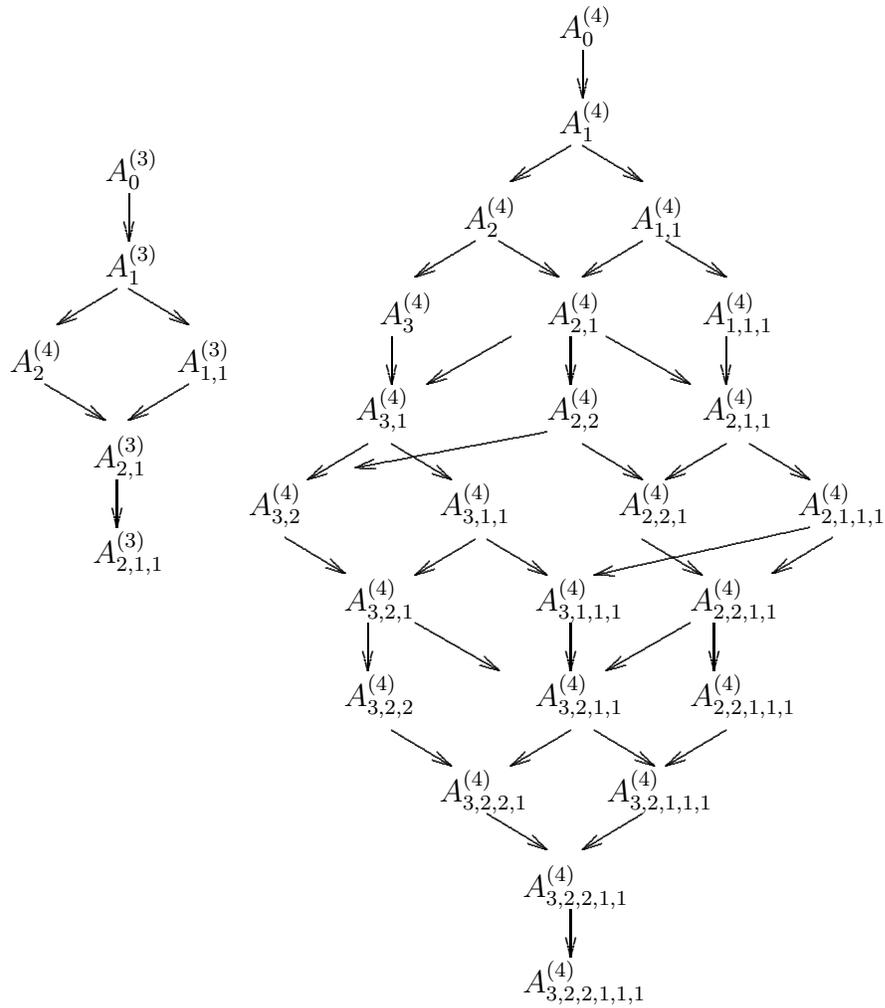}
\end{center}
\caption{Action of $e_1$ on irreducible atoms of level 3 and 4}
\label{figindecomp}
\end{figure}

We now have that
the dimension of the quotient $V_k(1)/\mathcal I_k$ is $k!$.  
Since we assume that the atoms of level $k$ form
a basis for $V_k$,  Corollary~\ref{conjrecschur} 
implies that the $k$-irreducible atoms 
also form a basis of $V_k(1)/\mathcal I_k$,
since the atoms generate $V_k(1)/\mathcal I_k$,  and the only possibly non-zero atoms in 
$V_k(1)/\mathcal I_k$ are the $k!$ irreducible ones.
\begin{corollary}  Assuming the atoms of level $k$ form a basis of $V_k$ and
Conjecture~\ref{conjrec} holds, the $k$-irreducible atoms 
form a basis of the quotient ring $V_k(1)/\mathcal I_k$.
\end{corollary}
If we link all atoms that occur in the action of $e_1$ on a given atom 
in $V_k(1)/\mathcal I_k$, we obtain a poset
illustrated in Figure~\ref{figindecomp}.
The rank generating function of this poset 
was given in \ref{HS}.  
This poset seems to have a remarkable 
symmetry property called {\it flip}-invariance.
\begin{definition}  
Given a $k$-irreducible partition of the form
\begin{equation}
\lambda = \bigl( (k-1)^{n_1},(k-2)^{n_2},\cdots,1^{n_{k-1}} \bigr) 
\quad\text{with}\;\; n_i \leq i\;\;\text{for all }\; i\, ,
\end{equation}
the involution called flip $\bf {f^{(k)}}$ 
is defined by
\begin{equation}
{\bf {f^{(k)}}} A_{\lambda}^{(k)} = A_{\lambda^{{\bf {f^{(k)}}}}}^{(k)} \,
\end{equation}
\begin{equation}
\text{where}\quad
\qquad \qquad
\lambda^{{\bf {f^{(k)}}}} = \bigl( (k-1)^{1-n_1},(k-2)^{2-n_2},\cdots,1^{k-1-n_{k-1}} \bigr) \, .
\qquad \qquad
\qquad \qquad
\qquad \qquad
\end{equation}
\end{definition}
\noindent For instance, 
\begin{equation*}
{\bf {f^{(5)}}} A_{4,3,2}^{(5)} = A_{3,2,2,1,1,1,1}^{(5)} \, .
\end{equation*}
\begin{conjecture}
The poset associated to the action of $e_1$ on atoms in $V_k(1)/\mathcal I_k$ is 
flip-invariant. 
That is, if there is an arrow between two atoms $A_{\mu}^{(k)}$ 
and $A_{\lambda}^{(k)}$, then there will be an arrow between 
the two atoms $A_{\mu^{\bf {f^{(k)}}}}^{(k)}$ 
and $A_{\lambda^{\bf {f^{(k)}}}}^{(k)}$.
\end{conjecture}

Given the $k!$ irreducible atoms, from
which all other atoms are constructed
using $k$-rectangular promotion operators,
the complete decomposition of the standard
tableaux into atoms can in principle be obtained.
We give here the cases $k=2$ and $k=3$. 

\subsection{Case $k=2$ and $k=3$}

We start with $k=2$. 
If $\mathcal S_n$ denotes the set of standard tableaux on $n$ letters,
then
\begin{equation}
\left(
\overline{\mathbb B}_{(2)} +
\overline{\mathbb B}_{(1^2)}
 \, \, \right)
\mathcal S_n = \mathcal S_{n+2} \, ,
\end{equation}
where
$\overline{\mathbb B}_{(2)}=\mathbb L_{\tiny{\tableau*[scY]{  1 & 2}}}\,
\mathbb B_{(2)}$ and $\overline{\mathbb B}_{(1^2)}=\mathbb B_{(1^2)}$.
This recursion implies,
for $\mathbb A_0^{(2)}=\mathbb H_0$ and
$\mathbb A_1^{(2)}=\tiny{\tableau*[scY]{1 \cr }}\,$,
\begin{equation}
\left(
\overline{\mathbb B}_{(2)} +
\overline{\mathbb B}_{(1^2)}
 \right)^\ell
\mathbb A_{\epsilon}^{(2)} = \mathcal S_{2\ell+\epsilon} \, ,
\quad\text{ where} \; \epsilon\in\{0,1\}\,.
\end{equation}
Expanding the left hand side gives
\begin{equation}
\sum_{(v_1,\dots,v_m)}
\,
\overline {\mathbb B}_{(v_1^{3-v_1})} \cdots
\overline {\mathbb B}_{ (v_m^{3-v_m})} \, \mathbb A_{\epsilon}^{(2)} \,
=\, \mathcal S_{2\ell+\epsilon} \, ,
\quad\text{for} \; v_i \in\{1,2\}\,,
\end{equation}
and each of the standard tableaux must occur in
exactly one term of this sum.
This is, each standard tableau must occur in exactly one family,
denoted
\begin{equation}
\mathbb A_{(v_1,\dots,v_m,\epsilon)}^{(2)} =
\overline {\mathbb B}_{(v_1^{3-v_1})} \cdots
\overline {\mathbb B}_{ (v_m^{3-v_m})} \, \mathbb A_{\epsilon}^{(2)}\, , \quad \quad  v_i \in\{1,2\}\,.
\end{equation}
We have thus decomposed the set of standard tableaux
into these families, which are the atoms of level 2
by Conjecture~\ref{conjrec}.  Furthermore, from Conjecture~\ref{conjkata},
given a standard tableau, we can determine to which family
$\mathbb A_{(v_1,\dots,v_m,\epsilon)}^{(2)}$ belongs, by first 
performing a (2)-katabolism ($v_1=2$)
if it contains the subword (12) and otherwise a (1,1)-katabolism ($v_1=1$).
Repeating this procedure on the resulting tableau (until there 
is only one box left ($\epsilon=1$) or no boxes left ($\epsilon=0$)),
we obtain the sequence $(v_1,\dots,v_m)$ that we need.

Now, from Conjecture~\ref{conjrec}, 
\begin{equation}
\digamma \left( \mathbb A_{(v_1,\dots,v_m,\epsilon)}^{(2)} \right) =
t^* A_{\lambda}^{(2)}[X;t] \, ,
\end{equation}
where $\lambda$ is the partition rearrangement of
$(v_1^{3-v_1},\dots,v_m^{3-v_m},\epsilon)$ and $*$ is a power of $t$.  
The symmetric function analogues of ${\mathbb B}_{(2)}$
and ${\mathbb B}_{(1^2)}$ are the vertex operators  $B_{(2)}$ and
$B_{(1^2)}$ (see next subsection).  Therefore, Conjecture~\ref{conjrec}
suggests that 
\begin{equation}
 B_{(v_1^{3-v_1})} \cdots
 B_{ (v_m^{3-v_m})} \, A_{\epsilon}^{(2)}[X;t]
= t^* A_{\lambda}^{(2)}[X;t] \, , \qquad \epsilon = 0,1 \, ,
\end{equation}
where $\lambda$ is the partition rearrangement of
$(v_1^{3-v_1},\dots,v_m^{3-v_m},\epsilon)$ and $*$ is a power of $t$.
This conjecture connects the atoms to the Macdonald polynomials, since
the creation operators that build the Macdonald polynomials recursively
can be divided into the operators
$B_{(2)}$ and $B_{(1^2)}$ \cite{[LM],[Z]}.
The positive expansion of Macdonald polynomials
indexed by 2-bounded partitions (equivalently, partitions 
with $\ell(\lambda) \leq 2$) into atoms of level 2 is thus 
conjecturally the one given in \cite{[LM],[Z]}
(and to \cite{[S1]}, since \cite{[ZS]} proves the operators 
are related to the functions studied in \cite{[S1]}).

In the case $k=3$, we have the 8 irreducible atoms
 of 6 distinct shapes,
\begin{equation}
\begin{split}
& \mathbb A_0^{(3)} = \mathbb H_0 \, \, ; \quad \mathbb A_1^{(3)}= {\tiny{\tableau*[scY
]{ 1}}} \, \, ; \quad \mathbb A_{12}^{(3)}= {\tiny{\tableau*[scY]{ 1&2}}} \, \,; \quad
\mathbb A_{21}^{(3)}= {\tiny{\tableau*[scY]{2 \cr 1}}} \, \, ;\\
  \mathbb A_{312}^{(3)}= {\tiny{\tableau*[scY]{ 3 \cr 1&2}}} \, \,; & \quad  \mathbb A_
{213}^{(3)}= {\tiny{\tableau*[scY]{ 2 \cr 1&3}}} \, \,; \quad  \mathbb A_{4312}^{(3)}=
{\tiny{\tableau*[scY]{ 4 \cr 3 \cr 1&2}}} + {\tiny{\tableau*[scY]{  3 \cr 1&2&4}}} \, 
\,; \quad \mathbb A_{4213}^{(3)}= {\tiny{\tableau*[scY]{ 4 \cr 2 \cr 1&3}}} + {\tiny{\tableau*[scY]{  2 \cr 1&3&4}}} \, \, ,
\end{split}
\label{3irr}
\end{equation}
from which we can build any atom of
evaluation $(1,\ldots,1)$ using
the promotion operators:
\begin{equation}
\mathbb B_{ {\tiny{\tableau*[scY]{ 1 & 2 &3}}}} \, \, , \,
\mathbb B_{ {\tiny{\tableau*[scY]{ 1 & 2 &4}}}} \, \, , \,
\mathbb B_{ {\tiny{\tableau*[scY]{ 1 & 3 &4 }}}} \, \, , \,
\mathbb B_{ {\tiny{\tableau*[scY]{ 3  \cr 2 \cr1 }}}} \, \, , \,
\mathbb B_{ {\tiny{\tableau*[scY]{ 4  \cr 2 \cr1 }}}} \, \, , \,
\mathbb B_{ {\tiny{\tableau*[scY]{ 4  \cr 3 \cr1 }}}} \, \, , \,
\mathbb B_{ {\tiny{\tableau*[scY]{ 3 & 4 \cr 1 & 2}}}} \, \, ,
\mathbb B_{ {\tiny{\tableau*[scY]{ 2 &4 \cr 1 &3}}}} \, \, , \,
\mathbb B_{ {\tiny{\tableau*[scY]{ 3 & 5\cr 1 & 2}}}} \, \, , \,
\mathbb B_{ {\tiny{\tableau*[scY]{ 2 & 5 \cr 1 & 3}}}} \, \, .
\label{promoop}
\end{equation}
Here an operator indexed by a tableau $T$ of shape $R$ is 
$\mathbb L_T \mathbb B_R$ followed by the reindexation of the letters 
not in $T$ such that
the resulting tableaux are standard.  For instance,
\begin{equation}
\mathbb B_{ {\tiny{\tableau*[scY]{ 1 & 3 &4 }}}} \, \, \, {\tiny{\tableau*[scY]{ 2 \cr
1 }}} \, = \,  {\tiny{\tableau*[scY]{5 \cr 2 \cr 1 & 3 &4  }}} +  {\tiny{\tableau*[scY]
{ 2 \cr 1 & 3 &4 & 5}}} \, .
\end{equation}

Using \ref{3irr} and \ref{promoop},
we consider the sets of tableau
\begin{equation}
\mathbb A_{(T_1,\dots,T_m,T)}^{(3)} = \mathbb B_{T_1} \cdots \mathbb B_{T_m} \,
\mathbb A_T^{(3)} \, ,
\end{equation}
for sequences $(T_1,\dots,T_m,T)$ that obey the following rules
(read from right to left):
\begin{enumerate}
\item[{1.}]  ${\tiny{\tableau*[scY]{ 3 & 5 \cr 1 &2 }}}$ and ${\tiny{\tableau*[scY]{ 2
& 5 \cr 1 &3 }}}$ can only follow a tableau that contains the subtableau ${\tiny{\tableau*[
scY]{ 1  }}}$  .
\item[{2.}]  ${\tiny{\tableau*[scY]{ 1 & 2 &4 }}}$ and ${\tiny{\tableau*[scY]{ 1 & 3 &4
 }}}$ can only follow a tableau that contains the subtableau ${\tiny{\tableau*[scY]{ 2 \cr
1  }}}$  .
\item[{3.}]  ${\tiny{\tableau*[scY]{ 4 \cr 2 \cr 1 }}}$ and ${\tiny{\tableau*[scY]{ 4 
\cr 3 \cr 1 }}}$ can only follow a tableau that contains the subtableau ${\tiny{\tableau*[sc
Y]{ 1 & 2  }}}$  .
\end{enumerate}
We conjecture that there is a one-to-one correspondence
between the sequences $(T_1,\dots,T_m,T)$,
and the set of tableaux indexing all level 3 copy atoms
with standard evaluation.  Moreover, we can determine to
which atom an arbitrary standard tableaux belongs
in view of Conjecture~\ref{conjkata};
katabolism is the inverse of rectangular promotion.
That is, given a tableau $U$,
we can determine which sequence $(T_1,\dots,T_m,T)$
can be extracted by katabolism from $U$.

\subsection{Generalized Kostka polynomials}

\label{secrec}

Given a sequence of partitions
$S=(\lambda^{(1)},\lambda^{(2)},\dots,\lambda^{(m)})$, 
the generalized Kostka polynomial $ H_{S}[X;t]$ found in \cite{[S1]}
is a $t$-generalization of the product of Schur functions 
indexed by the partitions in $S$ 
(different approaches to these polynomials include those in \cite{[LLT],[Shi]}).
More precisely, if we consider only its term of degree $n=|\lambda^{(1)}|+|\lambda^{(2)}|+\cdots+|\lambda^{(m)}|$,
\begin{equation}
 H_{S}[X;t] = \sum_{\lambda\vdash n} K_{\lambda;S}(t) \, S_{\lambda}[X] \, ,
\label{poinc}
\end{equation}
where, for the scalar product $\langle \, , \, \rangle$ on which the 
Schur functions are orthonormal,
\begin{equation}
K_{\lambda;S}(1) = 
\langle S_{\lambda}[X],S_{\lambda^{(1)}}[X] \, S_{\lambda^{(2)}}[X] \cdots 
\rangle \, .
\end{equation}
If successively reading the entries of 
$\lambda^{(1)}$, $\lambda^{(2)},\ldots$ produces a partition $\mu$, 
$S$ is said to be dominant.
In this case, it has been conjectured \cite{[S1]} that
\begin{equation}
 H_{S}[X;t]= \sum_{T \in {\mathbb H}_S } t^{\charge (T)} \, S_{\shape (T)}[X]\, ,
\end{equation}
where ${\mathbb H}_S$ is the set of tableaux $T$ of evaluation $\mu$ 
such that $\mathbb P_S(T)=T$ (see section~\ref{secatom}).
Now, if $S=\bigl((\ell_1^{k+1-\ell_1}), \dots,(\ell_m^{k+1-\ell_m}) \bigr)$ 
is a dominant sequence of $k$-rectangles, 
then Conjecture~\ref{conjrec} implies that for 
$\mu =\bigl(\ell_1^{k+1-\ell_1}, \dots,\ell_m^{k+1-\ell_m} \bigr)$,
\begin{equation}
\mathbb B_{(\ell_1^{k+1-\ell_1})} \cdots \mathbb B_{(\ell_m^{k+1-\ell_m})} \, \mathbb H_0  = \mathbb A^{(k)}_{\mu} \, .
\label{prodrec}
\end{equation}
Moreover, by the definition of atoms we have that
$\mathbb P_S ( \mathbb A^{(k)}_{\mu}) = \mathbb A^{(k)}_{\mu}$
since $\mu^{\to k}=S$.
Therefore, $\mathbb H_S = \mathbb A^{(k)}_{\mu}$
since both sets contain the same number of elements 
(the number of terms in the product of the Schur functions 
corresponding to shapes $(\ell_1^{k+1-\ell_1}), \dots,(\ell_m^{k+1-\ell_m})$).
We thus have the following connection between atoms
and the generalized Kostka polynomials:
\begin{conjecture}  If $S=\bigl((\ell_1^{k+1-\ell_1}), \dots,(\ell_m^{k+1-\ell_m}) \bigr)$ is such that $\bigl(\ell_1^{k+1-\ell_1}, \dots,\ell_m^{k+1-\ell_m} \bigr)$ is a partition $\mu$,  then 
\begin{equation}
A_{\mu}^{(k)}[X;t] = H_{S}[X;t] \, .
\end{equation}
\end{conjecture}

Further, it is shown in \cite{[ZS]} 
that the generalized Kostka polynomials can be defined as
\begin{equation}
 H_{S}[X;t] =   B_{\lambda^{(1)}}   B_{\lambda^{(2)}} \cdots  B_{\lambda^{(m)}} \, \cdot 1 \, ,
\label{poincop}
\end{equation}
where $B_{\lambda}$ corresponds to $H_{\lambda}^t$ in their notation.
Given our formula \ref{prodrec}, 
it is natural to assume that the vertex operators
$B_{(\ell^{k+1-\ell})}$ indexed by $k$-rectangular partitions 
are the operators that extend Conjecture~\ref{conjrec} to the 
level of symmetric functions.
\begin{conjecture}  
Given a  $k$-rectangular partition $(\ell^{k+1-\ell})$, we have
\begin{equation}
B_{(\ell^{k+1-\ell})} \, A_{\lambda}^{(k)}[X;t] = 
t^c \, A_{\lambda\cup(\ell^{k+1-\ell})}^{(k)}[X;t] \,,\quad
\text{where}\;\;c\in \mathbb N \,.
\end{equation}
\end{conjecture}

\section{The $k$-conjugation of a partition}

Here we introduce a generalization of partition conjugation,
defined for partitions bounded by $k$.  When $k$ is large,
our $k$-conjugation reduces to the usual conjugation.

A skew diagram $D$ is said to have hook-lengths bounded by $k$ if
the hook-length of any cell in $D$ is not larger than $k$. 
For a positive integer $m\leq k$, the
$k$-multiplication $m \times^{(k)} D$
is the skew diagram $\overline D$ obtained by adding 
a first column of length $m$ to $D$ such that the number of parts
of $\overline D$ is as small as possible while ensuring 
that its hook-lengths are bounded by $k$.
For example, 
\begin{equation}
 {\tiny{\tableau*[scY]{ \cr \cr \cr \cr}}} \, \times^{(5)} \,  {\tiny{\tableau*[scY]{ \cr \cr & \cr \bl & \cr \bl & & & \cr \bl & \bl & \bl & & }}} \, = \,  {\tiny{\tableau*[scY]{ \cr \cr & \cr & \cr\bl & & \cr\bl & \bl & \cr \bl & \bl & & & \cr\bl & \bl & \bl & \bl & & }}} \, .
\end{equation}

\begin{definition}
Let $\lambda=(\lambda_1,\dots,\lambda_n)$ be a $k$-bounded partition
and let $D$ be the skew diagram obtained by $k$-multiplying from right to 
left the entries of $\lambda$:
\begin{equation}
D =  \lambda_1 \times^{(k)} \cdots  \times^{(k)}
 \lambda_n  \, .
\label{spmul}
\end{equation}
The $k$-conjugate of $\lambda$,
denoted $\lambda^{\omega_k}$, 
is the vector obtained by reading the number of boxes in each row
of $D$.
\label{defconju}
\end{definition}

\begin{figure}[htb]
\begin{center}
\epsfig{file=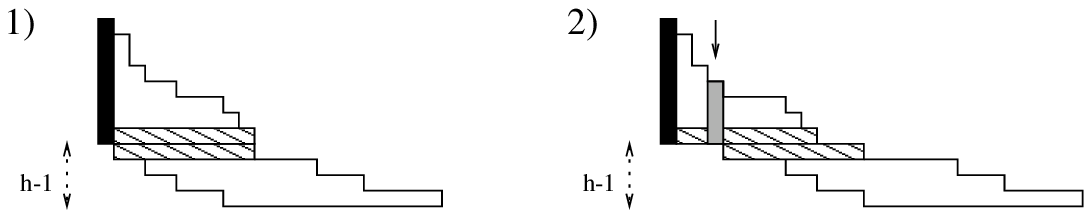}
\end{center}
\caption{}
\label{figconju}
\end{figure}

When $k \to \infty$, $\lambda^{\omega_k}=\lambda'$ 
since each $k$-multiplication step reduces to 
adding a column of length $\lambda_i$ at the bottom row.
\begin{property}
If $\lambda$ is a $k$-bounded partition,
then $\lambda^{\omega_k}$ is also a $k$-bounded partition.
\end{property}
\noindent {\it Proof.} \quad
$\lambda^{\omega_k}$ is $k$-bounded since
$D$ has hook-lengths bounded by $k$.  To see that 
$\lambda^{\omega_k}$ 
is a partition, assume by induction that the parts of
$D^{(2)} =  \lambda_2 \times^{(k)} \cdots  \times^{(k)} \lambda_n$ 
form a partition $\mu$. 
The skew diagram $D=\lambda_1 \times^{(k)} D^{(2)}$ is obtained
by adding a column of length $\lambda_1$ to $D^{(2)}$ starting at 
some row $h$.  
To see that $D$ must also have parts of weakly decreasing size,
it suffices to show that $\mu_{h-1}>\mu_{h}$.
Suppose $\mu_{h-1}=\mu_{h}$ and consider the two possible cases
(Figure~\ref{figconju}). Keep in mind that any 
column can be no longer than those to its left since 
$\lambda_i\leq\lambda_j\,\forall i>j$. 
If row $h-1$ lies directly below row $h$, then sliding 
the new column down to row $h-1$ 
gives a skew diagram of length less than $D$
with hook-lengths at most $k$.  Therefore 
our column would not have been added to row $h$.
Now if row $h-1$ lies below and to the right of row $h$,
the column indicated by an arrow can be moved down 
without producing any hook-lengths longer than $k$.  
Since $D^{(2)}=\lambda_2 \times^{(k)} \cdots  \times^{(k)} \lambda_n$, 
this is a contradiction.\hfill$\endprf$   

For example, 
we can compute $(2,2,1,1)^{\omega_4}=(3,2,1)$
by the following steps: 
\begin{equation}
 {\tiny{\tableau*[scY]{ \cr \cr }}} \, \times^{(4)} \,  
{\tiny{\tableau*[scY]{ \cr \cr  }}} \, \times^{(4)} \,  
{\tiny{\tableau*[scY]{ \cr}}} \, \times^{(4)} \,  
{\tiny{\tableau*[scY]{ \cr }}} \, = \, 
{\tiny{\tableau*[scY]{ \cr \cr }}} \, \times^{(4)} \,  
{\tiny{\tableau*[scY]{ \cr \cr  }}} \, \times^{(4)} \,  
{\tiny{\tableau*[scY]{ & }}} \, = 
\, \, {\tiny{\tableau*[scY]{ \cr \cr }}} \, \times^{(4)} \,  
{\tiny{\tableau*[scY]{  \cr & & \cr  }}} \, = 
\, {\tiny{\tableau*[scY]{  \cr &  \cr \bl &  & & \cr   }}} \, .
\label{spmulex}
\end{equation}

\begin{figure}[htb]
\begin{center}
\epsfig{file=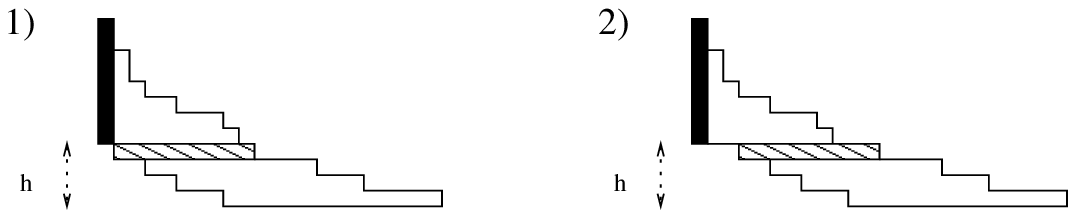}
\end{center}
\caption{}
\label{figinvo}
\end{figure}

\begin{property}  
For a $k$-bounded partition $\lambda$,
let $D =  \lambda_1 \times^{(k)} \cdots  \times^{(k)} \lambda_n$
and $\overline D$ be the skew diagram obtained by shifting any row
in $D$ to the left.
If the number of columns of $\overline D$ is 
not more than the number of columns of $D$ then 
the hook-lengths of $\overline D$ are not $k$-bounded.
\label{propmove}
\end{property}
\noindent {\it Proof.} \quad Assume by induction that 
$D^{(2)} =  \lambda_2 \times^{(k)} \cdots  \times^{(k)} \lambda_n$, 
with rows of length $\mu$, satisfies this 
property.  The skew diagram $D=\lambda_1 \times^{(k)} D^{(2)}$ falls into one
of the two generic cases illustrated in Figure~\ref{figinvo}.
In the first case, since the column is added above
row $h$,  we know $\lambda_1+\mu_h>k$.
Thus, row $h$ cannot be moved left
or we would have a cell with hook-length $\lambda_1+\mu_h>k$.  
In the second case, row $h$ cannot be moved left without 
violating the assumption that $D^{(2)}$ obeys the property. \hfill $\endprf$

\begin{theorem}  
$\omega_k$ is an involution on partitions bounded by $k$.  
That is, for $\lambda$ with $\lambda_1\leq k$,
\begin{equation}
\left( \lambda^{\omega_k} \right)^{\omega_k} = \lambda \, .
\end{equation}
\end{theorem}
\noindent {\it Proof.} \quad  
Let $D =  \lambda_1 \times^{(k)} \cdots  \times^{(k)} \lambda_n$.  
Property~\ref{propmove} implies that $D$ is recovered by
performing the $k$-multiplication of the entries of $\lambda^{\omega_k}$ 
in a conjugate way (adding rows to the leftmost position 
such that the hook-lengths are never larger than $k$).
Therefore, if $\lambda^{\omega_k}=\mu$, 
the conjugate of $D$ is given by $D'= \mu_1 \times^{(k)} \cdots  \times^{(k)} \mu_m$,
and thus $(D')'=D$ implies that 
$\mu^{\omega_k}= \left( \lambda^{\omega_k} \right)^{\omega_k}=\lambda$. 
\hfill $\endprf$

Given the $k$-conjugation of a partition, 
it is natural to consider
the relation among an atom indexed
by $\lambda$ and the atom indexed by
$\lambda^{\omega_k}$.
In fact, our examples suggest that
conjugating each tableaux in an
atom produces the tableaux
in another atom.
\begin{conjecture} Let $T$ be a standard tableau. 
For any copy $\mathbb A_T^{(k)}$ of 
$\mathbb A_{\lambda}^{(k)}$,
\begin{equation}
 \bigl( \mathbb A_T^{(k)} \bigr)^t  =  \mathbb A_{T'}^{(k)} \, ,
\end{equation}
for some standard tableau $T'$ of shape $\lambda^{\omega_k}$.
\label{conjtrans}
\end{conjecture}
Since, at any level, there is at least one copy of each atom of a given degree 
in the set of standard tableaux, we have the following corollary:
\begin{corollary} 
In any atom of shape $\lambda$ and level $k$, 
there is a unique element of maximal charge
whose shape is 
the conjugate of ${\lambda^{\omega_k}}$, 
$\left({\lambda^{\omega_k}}\right)'$. 
\label{cormin}
\end{corollary}

Furthermore, since a standard tableau $T$ in $n$ letters
satisfies $\charge(T^t)\!=\!{ n \choose 2}\! -\!\charge(T)$, 
\begin{corollary} 
Let $\omega$ be the involution such that $\omega S_{\lambda}[X]=S_{\lambda'}[X]$. 
Then, for some $*\in\mathbb N$,
\begin{equation}
\omega \, A_{\lambda}^{(k)}[X;t] =  t^*  A_{\lambda^{\omega_k}}^{(k)}[X;1/t] \, .
\end{equation}
\end{corollary}
\noindent Here we see that for large $k$,
$\lambda^{\omega_k}=\lambda'$
is consistent with the fact that 
$A_{\lambda}^{(k)}[X;t]=S_{\lambda}[X]$
in this case.

\section{Pieri rules}

Beautiful combinatorial algorithms are known for the
Littlewood-Richardson coefficients that
appear in a product of Schur functions;
\begin{equation}
S_{\lambda} \, S_{\mu} = \sum_{\nu} c_{\lambda \mu}^{\nu} \, S_{\nu} \, .
\label{schurpieri}
\end{equation}
Recall by Property~\ref{corlarge}
that our atoms $A_\lambda^{(k)}[X;t]$ are simply the
Schur functions $S_\lambda$ when $k$ is large.
Therefore the expansion coefficients in a product of
atoms are the Littlewood-Richardson coefficients when
$k$ is large and it is natural to examine the coefficients
in a product of two atoms for general $k$.
In fact, in the case $t=1$,
the coefficients in a product of two atoms do
seem to generalize Littlewood-Richardson coefficients.
\begin{conjecture}  Let $A_\lambda^{(k)}$ denote the
case $t=1$ in $A_\lambda^{(k)}[X;t]$.  Then
\begin{equation}
A_{\lambda}^{(k)} \, A_{\mu}^{(k)} =
\sum_{\nu} {c_{\lambda \mu}^{\nu}}^{\!\!(k)} A_{\nu}^{(k)} \, ,
\quad\text{where}\;\; 0 \subseteq
{c_{\lambda \mu}^{\nu}}^{\!\!(k)} \subseteq c_{\lambda \mu}^{\nu}
\;.
\end{equation}
In particular, we know
\begin{equation}
{c_{\lambda \mu}^{\nu}}^{\!\!(k)} = {c_{\lambda \mu}^{\nu}}
\quad
\text{for}\;\; k \geq |\mu| \, .
\end{equation}
\end{conjecture}

Identity \ref{schurpieri} reduces to the Pieri rule
when $\lambda$ is a row (resp. column).
Since an atom $A_{\lambda}^{(k)}$ reduces to $h_{\ell}$
(resp. $e_\ell$) when $\lambda$ is a row (column) of length $\ell \leq k$,
our conjecture can be reduced to a $k$-generalization of the Pieri rule.
\begin{corollary}
For certain sets of shapes $E_{\lambda,\ell}^{(k)}$ and $\bar E_{\lambda,\ell}^{(k)}$, we have for $\ell \leq k$,
\begin{equation}
h_\ell\, A_{\lambda}^{(k)} = \sum_{\mu \in E_{\lambda,\ell}^{(k)}}
A_{\mu}^{(k)}
\quad\text{and}
\quad
e_\ell\, A_{\lambda}^{(k)} = \sum_{\mu \in \bar E_{\lambda,\ell}^{(k)}}
A_{\mu}^{(k)} \, .
\end{equation}
\end{corollary}

We conjecture the sets 
$E_{\lambda,\ell}^{(k)}$ and $\bar E_{\lambda,\ell}^{(k)}$
can be defined in a manner analogous to the Pieri rule.
\begin{conjecture}
For any positive integer $\ell \leq k$,
\begin{equation}
\begin{split}
E_{\lambda,\ell}^{(k)}\, & = \,
\left\{\mu \, | \,
\mu/\lambda\; is \; a \; horizontal \;\ell\text{-}strip
\;\; and \;\;
\mu^{\omega_k}/\lambda^{\omega_k} \; is \; a \;  vertical
\; \ell\text{-}strip \right\} \, , \\
\bar E_{\lambda,\ell}^{(k)}\, & = \,
\left\{\mu \, | \, 
\mu/\lambda\; is \; a \; vertical \;\ell\text{-}strip
\;\; and \;\;
\mu^{\omega_k}/\lambda^{\omega_k} \; is \; a \;  horizontal
\; \ell\text{-}strip \right\} \, .
\end{split}
\end{equation}
\label{conjpieri}
\end{conjecture}
For example, to obtain the indices of the elements
that occur in $e_2 \, A_{3,2,1}^{(4)}$,
we compute $(3,2,1)^{\omega_4}=(2,2,1,1)$
by Definition~\ref{defconju}
and then add a horizontal 2-strip to (2,2,1,1)
in all possible ways.  This gives
(2,2,2,1,1),(3,2,1,1,1),(3,2,2,1) and (4,2,1,1)
of which all are 4-bounded.
Our set then consists of all the $4$-conjugates of these partitions
that leave a vertical 2-strip when 
$(3,2,1)$ is extracted from them. The corresponding 4-conjugates are
\begin{equation}
(2,2,2,1,1)^{\omega_4}   =  \,
{\tiny{\tableau*[scY]{ & \cr & & \cr & & \cr }}} \, ,
\,
(3,2,1,1,1)^{\omega_4}
= \, {\tiny{\tableau*[scY]{ \cr & \cr & \cr & & \cr }}}\, ,
\,
(3,2,2,1)^{\omega_4} = \, {\tiny{\tableau*[scY]{ \cr \cr \cr & \cr & & \cr }}} \, ,
\,
(4,2,1,1)^{\omega_4}  =  \, {\tiny{\tableau*[scY]{ \cr \cr \cr \cr \cr & & \cr }}} \,,
\end{equation}
and of these partitions,
only the first three are such that a
vertical 2-strip remains when $(3,2,1)$ is extracted.
Therefore
\begin{equation}
e_2 \, A^{(4)}_{3,2,1} \, = \,   A^{(4)}_{3,3,2} +
A^{(4)}_{3,2,2,1}  +  A^{(4)}_{3,2,1,1,1} \, ,
\end{equation}
which is in fact correct.

\section{Hook case}

We are able to explicitly determine the functions
$A_\lambda^{(k)}[X;t]$ in the case that $\lambda$ is
a hook partition and also to derive properties of
atoms indexed by partitions slightly more general
than hooks.  These results rely
on the following property of a row-shaped katabolism.
\begin{property}
If $T$ has shape $\lambda=(m,1^r)$ (a hook),
then
\begin{equation}
\mathbb K_{(n)} \, : \, T \longrightarrow
\begin{cases}
\bar T & {\text{if}\quad n\leq m}\\
0 & {\text{otherwise}} \\
\end{cases} \, \, ,
\end{equation}
where $\bar T$ is also hook-shaped.
\label{prophook}
\end{property}
\noindent {\it Proof.} \quad
Consider a tableau $T$ of shape $\lambda=(m,1^r)$.
If $n>m$ then $T$ does not contain a row of length $n$
and thus $\mathbb K_{(n)} T=0$.
Assume $n\leq m$. Let $U$ be the tableau of shape $(1^r)$
obtained by deleting the bottom row of $T$.
By the definition of katabolism,
the action of $\mathbb K_{(n)}$ on $T$
amounts to row inserting a sequence of
strictly decreasing letters (those of $U$)
into a sequence of weakly increasing letters
(the last $m-n$ letters in the bottom row of $T$).
The insertion algorithm implies \cite{[Tab]} that in this case,
no two elements may be added to the same row
and therefore, we obtain a hook shape.
\hfill $\endprf$

This property leads to the hook content of any
atom that is not indexed by a $k$-generalized hook
partition, that is, a partition
of the form $(k,\dots,k,\rho_1,\rho_2,\ldots)$
for a hook shape $(\rho_1,\rho_2,\ldots)$.
\begin{property}
If $T$ is a tableau of shape $\lambda$,
where $\lambda$ is not a $k$-generalized hook,
then $\mathbb A_T^{(k)}$
does not contain any tableaux with a hook shape.
\end{property}
\begin{corollary}
If $\lambda$ is a partition that is not a $k$-generalized hook,
then
\begin{equation}
A_{\lambda}^{(k)}[X;t] = S_{\lambda}[X] +
\sum_{\mu > \lambda} v_{\mu \lambda}^{(k)}(t) \, S_{\mu}[X]\, ,
\end{equation}
where $v_{\mu\lambda}^{(k)}(t)=0$ for all hook partitions $\mu$.
\label{atomtri}
\end{corollary}
\noindent {\it Proof.} \quad
Let $\lambda^{\to k}=(\lambda^{(1)},\lambda^{(2)},\ldots)$.
The condition on $\lambda$ implies that $\lambda_2$ is at least 2.
If we first consider such partitions with $\lambda_1\neq k$
then $\lambda^{(1)}$ cannot be a hook
($\lambda_1\neq k$ implies that the first partition in the $k$-split
contains at least the first two parts of $\lambda$).
But if $\lambda^{(1)}$ is not a hook,
then any hook-shaped tableau $T$ in
$\mathbb B_{(\lambda_1)}\left( \mathbb A_{(\lambda_2,\lambda_3,\dots)}^{(k)}
\right)$  will not contain the shape $\lambda^{(1)}$
and will therefore be sent to zero under $\mathbb P_{\lambda^{\to k}}$.
On the other hand, if $\lambda_1=k$ then $\lambda^{(1)}=(k)$.
Now any hook-shaped tableau $T$ in
$ \mathbb B_{(\lambda_1)}\left( \mathbb A_{(\lambda_2,\lambda_3,\dots)}^{(k)}
\right)$ will be sent to a hook under the $(k)$-katabolism by
Property~\ref{prophook}. Our claim thus
follows recursively on the remaining terms of the $k$-split
of $\lambda$.
\hfill $\endprf$

If an atom {\it is} indexed by a $k$-generalized hook, we can determine it explicitly.
\begin{property}
Let $\lambda=(m,1^r)$ be a $k$-irreducible hook partition.  Then
\begin{equation}
\mathbb A_{\lambda}^{(k)} =
\begin{cases}
(r+1) \, r \cdots 2 \, 1^m & {\rm if~} r +m \leq k \\
 (r+1) \, r \cdots 2 \, 1^m +  r \cdots 2 \, 1^m \, (r+1) & {\rm otherwise}
\end{cases} \, .
\end{equation}
Note, here an element $(r+1)\,r\cdots 2\,1^m$ denotes the word
$(r+1)\,r\cdots 2\,1\,1\cdots 1$.
\end{property}

\noindent {\it Proof.} \quad
Since $r,m \leq k-1$ in any $k$-irreducible partition $\lambda=(m,1^r)$,
we have that $(1^i)^{\to k}=(1^i)$  for $1\leq i\leq r$.
Therefore, on a tableau $T$ with $i$ boxes,
$\mathbb P_{(1^i)^{\to k}} T \neq 0 $ only for $T$ of shape $(1^i)$
and thus
\begin{equation}
\mathbb A_{1^r}=
\mathbb P_{(1^r)^{\to k}} \, \mathbb B_1
\cdots
\mathbb P_{(1^2)^{\to k}} \, \mathbb B_1
\mathbb P_{(1)^{\to k}} \, \mathbb B_1
\, \mathbb H_{0}
\, = \,
r\,\,r-1\cdots 1
\, .
\end{equation}
Moreover, it develops that
$\mathbb B_m \, \mathbb A_{1^r}= (r+1) \, r \cdots 2 \, 1^m \,
+ \, r \cdots 2 \, 1^m \, (r+1)$.
Now we have
\begin{equation}
\mathbb A_{(m,1^r)}=
\mathbb P_{(m,1^r)^{\to k}}
\big( (r+1) \, r \cdots 2 \, 1^m \,
+ \, r \cdots 2 \, 1^m \, (r+1)\big)\,.
\end{equation}
Since $r+m-k \leq k$,
the $k$-split of $(m,1^r)$ is
\begin{equation}
(m,1^r)^{\to k} =
\begin{cases}
\bigl( (m,1^r) \bigr) & {\rm if~} r+m \leq k \\
\bigl( (m,1^{k-m}), (1^{r+m-k}) \bigr) & {\rm otherwise}
\end{cases}
\, .
\end{equation}
In either case, $\mathbb P_{(m,1^r)^{\to k}}
\left( (r+1) \, r \cdots 2 \, 1^m\right)
= (r+1) \, r \cdots 2 \, 1^m$,
but since $r \cdots 2 \, 1^m \, (r+1)$
never contains shape $(m,1^r)$,
$\mathbb P_{(m,1^r)^{\to k}}
\left(r \cdots 2 \, 1^m \, (r+1)\right)\neq 0$
only in the second case.
\hfill $\endprf$

Now by Conjecture~\ref{conjrec},
we use the given atoms of level $k$ indexed by a $k$-irreducible hook shape
to obtain more general cases including those indexed by
a $k$-generalized hook shape.
\begin{corollary}  Assume Conjecture~\ref{conjrec} holds.
For a sequence of $k$-rectangles
$(R_1,R_2,\ldots,R_j)$,
let $\lambda$ be the partition rearrangement of
$(R_1,R_2,\ldots,R_j,m,1^r)$.
Then
\begin{equation}
A_{\lambda}^{(k)}[X;t] \propto
\digamma \left( 
\mathbb B_{R_1}
\cdots
\mathbb B_{R_j}
\mathbb A_{(m,1^r)} \right) 
\, .
\end{equation}
\end{corollary}

\end{document}